\documentclass[lang = american]{ems-icm}

\usepackage{amscd,mathtools, mathtext,color,wrapfig, color, verbatim, slashed}

\usepackage{tikz-cd}
\usepackage{bbm}
\usepackage{upgreek}
\usepackage{extarrows} 
\usepackage{diagbox}
\usepackage[mathscr]{eucal}
\usepackage{mathrsfs}
\usepackage{pifont}
\usepackage{cleveref}
\usepackage[font={small,sl}]{caption}
\usepackage[all]{xy}
\usepackage{wasysym}
\usepackage{tikz}
\usetikzlibrary{matrix,arrows}

\newcommand{\beq}{\begin{equation}}
\newcommand{\eeq}{\end{equation}}

\newcommand{\nc}{\newcommand}
\nc{\mc}{\mathcal}
\newtheorem{theorem}{Theorem}

\newtheorem{corrolary}{Corrolary}
\newenvironment{theorem*}
 {\expandafter\def\expandafter\thetheorem\expandafter{\thetheorem*}\theorem}
 {\endtheorem}
\newenvironment{theorem^!}
 {\expandafter\def\expandafter\thetheorem\expandafter{\thetheorem^!}\theorem}
 {\endtheorem}
\newenvironment{corrolary*}
 {\expandafter\def\expandafter\thecorrolary\expandafter{\thecorrolary*}\corrolary}
 {\endcorrolary}

\usetikzlibrary{arrows,chains,matrix,positioning,scopes}
\makeatletter
\tikzset{join/.code=\tikzset{after node path={%
\ifx\tikzchainprevious\pgfutil@empty\else(\tikzchainprevious)%
edge[every join]#1(\tikzchaincurrent)\fi}}}
\makeatother
\tikzset{>=stealth',every on chain/.append style={join},
         every join/.style={->}}
\tikzstyle{labeled}=[execute at begin node=$\scriptstyle,
   execute at end node=$]
\DeclareMathAlphabet{\pazocal}{OMS}{zplm}{m}{n}

\numberwithin{equation}{section}

\begin{document}
\volume{?}

%\baselineskip=28pt  % a la harvmac
%\baselineskip 0.7cm
%\baselineskip=18pt  % a la harvmac
%\baselineskip 0.7cm
%\setcounter{tocdepth}{2}

%\begin{titlepage}

%% Set the number of the title with 0

% change the footnote symbol
%\renewcommand{\thefootnote}{\fnsymbol{footnote}}

%\vskip 1.0cm
\title{Homological Knot Invariants from Mirror Symmetry}
\titlemark{}

%\begin{center}
%{\LARGE 
% Homological Knot Invariants 
%\vskip 0.5 cm
%from Mirror Symmetry
%\vskip 0.5cm}

\emsauthor{1}{Mina Aganagic}{M.~Aganagic}
%{\large
%Mina Aganagic}

%\medskip

%\vskip 0.5cm

{\emsaffil{1}{
Department of Mathematics, University of California, Berkeley\\
Center for Theoretical Physics, University of California, Berkeley
 \email{aganagic@berkeley.edu}
 }

\begin{abstract}
In 1999, Khovanov showed that a link invariant known as the Jones polynomial
is the Euler characteristic of a homology theory. The knot categorification
problem is to find a general construction of knot homology groups, and
to explain their meaning -- what are they homologies of?

Homological mirror symmetry, formulated by Kontsevich in 1994, naturally
produces hosts of homological invariants. Typically though, it leads to
invariants which have no particular interest outside of the problem at
hand.

I showed recently that there is a new family of mirror pairs of manifolds,
for which homological mirror symmetry does lead to interesting invariants
and solves the knot categorification problem. The resulting invariants
are computable explicitly for any simple Lie algebra, and certain Lie superalgebras.
\end{abstract}

%\classification[53D45, 53D37, 81]{14J33}
% one primary code (in curly braces)
        % and a list of secondary codes separated by commas (in square brackets)
%\keywords{Homological mirror symmetry, knot theory, categorification} %Only the first keyword is capitalized
%

\maketitle

\section{Introduction}
\label{197.sec1}

There are many beautiful strands that connect mathematics and physics.
Two of the most fruitful ones are knot theory and mirror symmetry.
I will describe a new connection between them. We will find a solution
to the knot categorification problem as a new application of homological
mirror symmetry.

%s1.1 #&#
\subsection{Quantum link invariants}
\label{197.sec1.1}

In 1984, Jones constructed a polynomial invariant of links in
${\mathbb R}^{3}$ \cite{JonesVFR}. The Jones polynomial is defined by picking
a projection of the link to a plane, the skein relation it satisfies
%
%\vspace*{9pt}
\begin{align*}
\includegraphics{197i01.pdf}
\end{align*}
where $n=2$, and the value for the unknot. It has the same flavor as the Alexander polynomial, dating
back to 1928 \cite{Alexander}, which one gets by setting $n=0$ instead.

The proper framework for these invariants was provided by Witten in 1988,
who showed that they originate from three-dimensional Chern--Simons theory
based on a Lie algebra $^{L}{\mathfrak{g}}$ \cite{Jones}. In particular,
the Jones polynomial comes from
$^{L}{\mathfrak{g}} = {\mathfrak{su}}_{2}$ with links colored by the defining
two-dimensional representation. The Alexander polynomial comes from the
same setting by taking $^{L}{\mathfrak{g}}$ to be a Lie superalgebra
${\mathfrak{gl}}_{1|1}$. The resulting link invariants are known as the
$U_{\mathfrak{q}}(^{L}{\mathfrak{g}})$ quantum group invariants. The relation
to quantum groups was discovered by Reshetikhin and Turaev \cite{RT}.

%s1.2 #&#
\subsection{The knot categorification problem}
\label{197.sec1.2}

The quantum invariants of links are Laurant polynomials in
${\mathfrak{q}}^{1/2}$, with integer coefficients. In 1999, Khovanov showed
\cite{Kh, KhICM} that one can associate to a projection of the link to
a plane a bigraded complex of vector spaces
\begin{equation*}
C^{*, j}(K) = \cdots C^{i-1, j} (K)\xrightarrow{\partial ^{i-1}} C^{i,j}(K)
\xrightarrow{\partial ^{i}} \cdots ,
\end{equation*}
whose cohomology
${\mathcal H}^{i, j}(K) = {\mathrm{ker} \,{\partial ^{i}}/
\mathrm{im}\, \partial ^{i-1}}$ categorifies the Jones polynomial,
\begin{equation*}
J_{K}({\mathfrak{q}}) = \sum _{i,j} (-1)^{i} {\mathfrak{q}}^{j/2}
{\mathrm{rk}} \,{\mathcal H}^{i,j}(K).
\end{equation*}
Moreover, the cohomology groups
\begin{equation*}
{\mathcal H}^{*,*}(K) = \bigoplus _{i,j} {\mathcal H}^{i,j}(K)%
\end{equation*}
are independent of the choice of projection; they are themselves link invariants.

%s1.2.1 #&#
\subsubsection{}
\label{197.sec1.2.1}

Khovanov's construction is part of the categorification program initiated
by Crane and Frenkel \cite{CF}, which aims to lift integers to vector spaces
and vector spaces to categories. A toy model of categorification comes
from a Riemannian manifold $M$, whose Euler characteristic
\begin{equation*}
{\chi (M)} = \sum _{k\in \mathbb{Z}} (-1)^{k} \mathrm{dim}\,{\mathcal H}^{k}(M)
\end{equation*}
is categorified by the cohomology
${\mathcal H}^{k}(M) = {\mathrm{ker} \;{d_{k}}/ \mathrm{im}\; d_{k-1}}$
of the de Rham complex
\begin{equation*}
C^{*} = \cdots C^{k-1} \xrightarrow{d_{k-1}} C^{k}
\xrightarrow{d_{k}} \cdots .%
\end{equation*}
The Euler characteristic is, from the physics perspective, the partition
function of supersymmetric quantum mechanics with $M$ as a target space
$\chi (M) = \mathrm{Tr} (-1)^{F} e^{-\beta H}$, with Laplacian $H=d d^{*}+ d^{*} d$ as the Hamiltonian, and $d=\sum _{k} d_{k}$ as the supersymmetry
operator. If $h$ is a Morse function on $M$, the complex can be replaced
by a Morse--Smale--Witten complex $C^{*}_{h}$ with the differential
$d_{h} = e^{h} d e^{-h}$. The complex $C^{*}_{h}$ is the space of perturbative ground states of a
$\sigma $-model on $M$ with potential $h$ \cite{WittenM}. The action of the differential
$d_{h} $ is generated by solutions to flow equations, called instantons.

%s1.2.2 #&#
\subsubsection{}
\label{197.sec1.2.2}

Khovanov's remarkable categorification of the Jones polynomial is explicit
and easily computable. It has generalizations of similar flavor for
$^{L}{\mathfrak{g}} = \mathfrak{su}_{n}$, and links colored by its minuscule
representations \cite{KR}.

In 2013, Webster showed \cite{webster} that for any
$^{L}{\mathfrak{g}}$, there exists an algebraic framework for categorification
of $U_{\mathfrak{q}}(^{L}{\mathfrak{g}})$ invariants of links in
${\mathbb R}^{3}$, based on a derived category of modules of an associative algebra.
The KLRW algebra, defined in \cite{webster}, generalizes the algebras of
Khovanov and Lauda \cite{KL_} and Rouquier \cite{R}. Unlike Khovanov's construction,
Webster's categorification is anything but explicit.

%s1.2.3 #&#
\subsubsection{}
\label{197.sec1.2.3}

Despite the successes of the program, one is missing a fundamental principle
which explains why is categorification possible -- the construction has
no right to exist. Unlike in our toy example of categorification of the
Euler characteristic of a Riemanniann manifold, Khovanov's construction
and its generalizations did not come from either geometry or physics in
any unified way. The problem Khovanov initiated is to find a general framework
for link homology, that works uniformly for all Lie algebras, explains
what link homology groups are, and why they exist.

%s1.3 #&#
\subsection{Homological invariants from mirror symmetry}
\label{197.sec1.3}

The solution to the problem comes from a new relation between mirror symmetry
and representation theory.

Homological mirror symmetry relates pairs of categories of geometric origin
\cite{Ko}: a derived category of coherent sheaves and a version of the
derived Fukaya category, in which complementary aspects of the theory are
simple to understand. Occasionally, one can make mirror symmetry manifest,
by showing that both categories are equivalent to a derived category of
modules of a single algebra.

I will describe a new family of mirror pairs, in which homological mirror
symmetry can be made manifest and leads to the solution to the knot categorification
problem \cite{A1, A2}. Many special features exist in this family, in part
due to its deep connections to representation theory. As a result, the theory is
solvable explicitly, as opposed to only formally \cite{ALR, ADZ}.

%s1.4 #&#
\subsection{The solution}
\label{197.sec1.4}

We will get not one, but two solutions to the knot categorification problem.
The first solution \cite{A1} is based on $\mathscr{D}_{\mathcal X}$, the
derived category of equivariant coherent sheaves on a certain holomorphic
symplectic manifold ${\mathcal X}$, which plays a role in the geometric
Langlands correspondence. Recently, Webster proved that
$\mathscr{D}_{\mathcal X}$ is equivalent to
$\mathscr{D}_{\mathscr A}$, the derived category of modules of an algebra
${\mathscr A}$ which is a cylindrical version of the KLRW algebra from
\cite{W1,W2}. The generalization allows the theory to describe links in
${\mathbb R}^{2}\times S^{1}$, as well as in ${\mathbb R}^{3}$.

The second solution \cite{A2} is based on $\mathscr{D}_{Y}$, the derived
Fukaya--Seidel category of a certain manifold $Y$ with potential $W$. The
theory generalizes Heegard--Floer theory \cite{OS0b, OS0c, Rasmussen}, which categorifies the Alexander
polynomial, from $^{L}{\mathfrak{g}} = {\mathfrak{gl}}_{1|1}$, to arbitrary
$^{L}{\mathfrak{g}}$.

The two solutions are related by equivariant homological mirror symmetry,
which is not an equivalence of categories, but a correspondence of objects
and morphisms coming from a pair of adjoint functors. In
$\mathscr{D}_{\mathcal X}$, we will learn which question we need to ask
to obtain $U_{\mathfrak{q}}(^{L}\mathfrak{g})$ link homology. In $\mathscr{D}_{Y}$, we will learn how to
answer it.

In \cite{ALR}, we give an explicit algorithm for computing homological link
invariants from $\mathscr{D}_{Y}$, for any simple Lie algebra
$^{L}{\mathfrak{g}}$ and links colored by its minuscule representations.
It has an extension to Lie superalgebras
$^{L}{\mathfrak{g}} = {\mathfrak{gl}}_{m|n}$ and
${\mathfrak{sp}}_{m|2n}$. In \cite{ADZ}, we set the mathematical foundations
of ${\mathscr{D}_{Y}}$ and prove (equivariant) homological mirror symmetry
relating it to ${\mathscr{D}_{\mathcal X}}$.

%s2 #&#
\section{Knot invariants and conformal field theory}
%%LEAP%%%\label{197.sec2}
\label{CFT}

Most approaches to categorification of
$U_{\mathfrak{q}}(^{L}{\mathfrak{g}})$ link invariants start with quantum
groups and their modules. We will start by recalling how quantum groups
came into the story. The seeming detour will help us understand how
$U_{\mathfrak{q}}(^{L}{\mathfrak{g}})$ link invariants arise from geometry,
and what categorifies them.

%s2.1 #&#
\subsection{Knizhnik--Zamolodchikov equation and quantum groups}
\label{197.sec2.1}

Chern--Simons theory associates to a punctured Riemann surface ${\mathcal A}$ a vector space, its Hilbert space. As Witten showed \cite{Jones}, the Hilbert
space is finite dimensional, and spanned by vectors that have a name. They are known as conformal blocks of the affine Lie
algebra $\widehat{^L\mathfrak{g}}_{\kappa}$. The effective level
$\kappa $ is an arbitrary complex number, related to
${\mathfrak{q}}$ by $ \mathfrak{q}= e^{\frac{2\pi i}{\kappa}}$. While in
principle arbitrary representations of $^{L}\mathfrak{g}$ can occur, in
relating to geometry and categorification, we will take them to be minuscule.

To get invariants of knots in ${\mathbb R}^{3}$, one typically takes
${\mathcal A}$ to be a complex plane with punctures. It is equivalent,
but for our purposes better, to take ${\mathcal A}$ to be an infinite
complex cylinder. This way, we will be able to describe invariants of links
in ${\mathbb R}^{2}\times S^{1}$, as well.

%s2.1.1 #&#
\subsubsection{}
\label{197.sec2.1.1}

Every conformal block, and hence every state in the Hilbert space, can
be obtained explicitly as a solution to a linear differential equation
discovered by Knizhnik and Zamolodchikov in 1984 \cite{KZ}. The KZ equation
we get is of trigonometric type, schematically
%
%e\upshape  2.1 #&#
\begin{equation}
\label{KZ}
\kappa \, {\partial _{i}} \,{\mathcal V} = \sum _{j \neq i} {r_{i j}}(a_{i}/a_{j})
\, {\mathcal V},
\end{equation}
since ${\mathcal A}$ is an infinite cylinder. Here,
$\partial _{i} = a_{i} {\frac{\partial}{\partial a_{i}}}$, where
$a_{i}$ is any of the $n$ punctures in the interior of
${\mathcal A}$, colored by a representation $V_{i}$ of
$^{L}{\mathfrak{g}}$. The right hand side of \eqref{KZ} is given in terms of classical
$r$-matrices of $^{L}{\mathfrak{g}}$, and acts irreducibly on a subspace of
$V_{1} \otimes \dots \otimes V_{n}$ of a fixed weight $\nu $, where
${\mathcal V}$ takes values \cite{EFK, EG}.

The KZ equations define a flat connection on a vector bundle over the
configuration space of distinct points
$a_{1}, \ldots , a_{n} \in {\mathcal A}$. The flatness of the connection
is the integrability condition for the equation.

%s2.1.2 #&#
\subsubsection{}
\label{197.sec2.1.2}

The monodromy problem of the KZ equation, which is to describe analytic
continuation of its fundamental solution along a path in the configuration space, has an explicit solution. Drinfeld
\cite{Drinfeld} and Kazhdan and Lustig \cite{KL} proved that that the monodromy matrix
${\mathfrak B}$ of the KZ connection is a product of $R$-matrices of the
$U_{\mathfrak{q}}(^{L}{\mathfrak{g}})$ quantum group corresponding to
$^{L}{\mathfrak{g}}$. The $R$-matrices describe reorderings of neighboring
pairs of punctures.

The monodromy matrix
${\mathfrak B}$ is the Chern--Simons path integral on
${\mathcal A} \times [0,1]$ in presence of a colored braid.
By composing braids, we get a representation of the affine braid group
based on the $U_{\mathfrak{q}}(^{L}{\mathfrak{g}})$ quantum group, acting on the space of solutions to the KZ equation. The
braid group is affine, since ${\mathcal A}$ is a cylinder and not a plane.

%s2.1.3 #&#
\subsubsection{}
\label{197.sec2.1.3}

%f\caption@xref {NNEEWWLLAABBEELL93(thesamelabel) 1 #&#
\begin{figure}[b!]
\begin{center}
\vspace*{-18pt}
\includegraphics{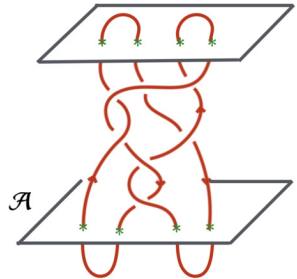}
\end{center}
\caption{Every link arises as a plat closure of a braid.}\label{f_1}
\end{figure}

Any link can be represented as a plat closure of some braid.  The Chern-Simons path integral together with the link computes a very specific matrix
element of the braiding matrix ${\mathfrak B}$, picked out by a pair of states in the Hilbert space corresponding to the top and the bottom of Figure \ref{f_1}.

These states, describing a collection of cups and caps, are very special solutions of the KZ equation in which
pairs of punctures, colored by conjugate representations $V_{i}$ and
$V_{i}^{*}$, come together and fuse to disappear. In this way, both
fusion and braiding enter the problem.

%s2.2 #&#
\subsection{A categorification wishlist}
\label{197.sec2.2}

To categorify $U_{\mathfrak{q}}(^{L}{\mathfrak{g}})$ invariants of links
in ${\mathbb R}^{3}$, we would like to associate, to the space of conformal
blocks of $\widehat{^L\mathfrak{g}}$ on the Riemann surface
${\mathcal A}$, a bigraded category, which in addition to the cohomological
grading has a grading associated to ${\mathfrak{q}}$.
Additional
$\mathrm{rk}(^{L}{\mathfrak{g}})$ gradings are needed to categorify invariants
of links in ${\mathbb R}^{2}\times S^{1}$, as they depend on the choice
of a flat $^{L}{\mathfrak{g}}$ connection around the $S^{1}$. To braids,
we would like to associate functors between the categories corresponding
to the top and bottom. To links, we would like to associate a vector space
whose elements are morphisms between the objects of the categories associated to the top
and bottom, up to the action of the braiding functor. Moreover, we would
like to do this in a way that recovers quantum link invariants upon decategorification.
One typically proceeds by coming up with a category, and then works to
prove that decategorification gives the link invariants one set out to
categorify. A virtue of the solutions in \cite{A1, A2} is that the second
step is automatic.

%s3 #&#
\section{Mirror symmetry}
\label{197.sec3}

Mirror symmetry is a string duality which relates $\sigma $-models on a pair
of Calabi--Yau manifolds ${\mathcal X}$ and ${\mathcal Y}$. Its mathematical
imprint are relations between very different problems in complex geometry
of ${\mathcal X}$ (``B-type'') and symplectic geometry of
${\mathcal Y}$ (``A-type''), and vice versa.

Mirror symmetry was discovered as a duality of $\sigma $-models on closed
Riemann surfaces ${D}$. In string theory, one must allow Riemann surfaces
with boundaries. This enriches the theory by introducing ``branes,'' which
are boundary conditions at $\partial{D}$ and naturally objects of a category
\cite{Douglas}.

By asking how mirror symmetry acts on branes turned out to yield deep insights
into mirror symmetry. One such insight is due to Strominger, Yau, and Zaslov
\cite{SYZ}, who showed that in order for every point-like brane on
${\mathcal X}$ to have a mirror on ${\mathcal Y}$, mirror manifolds have
to be fibered by a pair of (special Lagrangian) dual tori $T$ and
$T^{\vee}$, over a common base.

%s3.1 #&#
\subsection{Homological mirror symmetry}
\label{197.sec3.1}

Kontsevich conjectured in his 1994 ICM address \cite{Ko} that mirror symmetry
should be understood as an equivalence of a pair of categories of branes,
one associated to complex geometry of ${\mathcal X}$, the other to symplectic
geometry of ${\mathcal Y}$.

The category of branes associated to complex geometry of
${\mathcal X}$ is the derived category of coherent sheaves,
\begin{equation*}
{\mathscr D}_{\mathcal X} = D^b\mathrm{Coh}_{{T}}({\mathcal X}).
\end{equation*}
Its objects are ``B-type branes,'' supported on complex submanifolds of
${\mathcal X}$. The category of branes associated to symplectic geometry
is the derived Fukaya category
\begin{equation*}
{\mathscr D}_{\mathcal Y} =D\!\operatorname{Fuk}({\mathcal Y}),
\end{equation*}
whose objects are ``A-type branes,'' supported on Lagrangian submanifolds
of ${\mathcal Y}$, together with a choice of orientation and a flat bundle.
For example, mirror symmetry should map the structure sheaf of a point
in ${\mathcal X}$ to a Lagrangian brane in ${\mathcal Y}$ supported on
a $T^{\vee}$ fiber. The choice of a flat $U(1)$ connection is the position
of the point in the dual fiber $T$.

Kontsevich' homological mirror symmetry is a conjecture that the category
of B-branes on ${\mathcal X}$ and the category of A-branes on
${\mathcal Y}$ are equivalent,
\begin{equation*}
{\mathscr D}_{\mathcal X} \cong {\mathscr D}_{\mathcal Y},
\end{equation*}
and that this equivalence characterizes what mirror symmetry is.

%s3.2 #&#
\subsection{Quantum differential equation and its monodromy}
\label{197.sec3.2}

Knizhnik--Zamolodchikov equation, which plays a central role in knot theory,
has a geometric counterpart. In the world of mirror symmetry, there is
an equally fundamental differential equation,
%
%e\upshape  3.1 #&#
\begin{equation}
\label{qdif}
{\partial _{i}} {\mathcal V}_{\alpha} - (C_{i})_{\alpha}^{\beta} \,
{\mathcal V}_{\beta} =0.
\end{equation}
The equation is known as the quantum differential equation of
${\mathcal X}$. Both the equation and its monodromy problem featured prominently,
starting with the first papers on mirror symmetry, see \cite{G0} for an
early account. In \eqref{qdif},
$(C_{i})_{\alpha}^{\beta} = C_{\gamma _{i} \alpha \delta}
\eta ^{\delta \beta}$ is a connection on a vector bundle with fibers
$H^{\mathrm{even}}({\mathcal X}) = \bigoplus _{k} H^{k}({\mathcal X},\wedge ^{k}
T_{\mathcal X}^{*} )$ over the complexified Kahler moduli space.  The derivative stands for
$\partial _{i} = a_{i} {\frac{\partial}{\partial a_{i}}}$, so that
$\partial _{i} a^{d} = (\gamma _{i}, d) a^{d}$ for a curve of degree
$d\in H_{2}({\mathcal X})$. The connection comes
from quantum multiplication with classes
$\gamma _{i} \in {H}^{2}({\mathcal X})$. Given three de Rham cohomology classes on
${\mathcal X}$, their quantum product
%
%e\upshape  3.2 #&#
\begin{equation}
\label{qmulti}
C_{\alpha \beta \gamma} = \sum _{d\geq 0, d\in H_{2}({\mathcal X})} (
\alpha , \beta , \gamma )_{d} \; a^{d}
\end{equation}
is a deformation of the classical cup product
$ (\alpha , \beta , \gamma )_{0} =\int _{\mathcal X} \alpha \wedge
\beta \wedge \gamma $ coming from Gromov--Witten theory of
${\mathcal X}$: $ (\alpha , \beta , \gamma )_{d}$ is computed by an integral
over the moduli space of degree $d$ holomorphic maps from
${D} = {\mathbb P}^{1}$ to ${\mathcal X}$ whose image meets classes Poincar\'e
dual of $\alpha , \beta $ and $\gamma $ at points. The quantum product,
together with the invariant inner product
$\eta _{\alpha \beta} = \int _{\mathcal X} \alpha \wedge \beta $, gives
rise to an associative algebra with structure constants
${C_{\alpha \beta}}^{\delta} = C_{\alpha \beta \gamma} \eta ^{\gamma
\delta}$. Flatness of the connection follows from the WDVV equations
\cite{wt, DVV, Dubrovin}.

From the mirror perspective of ${\mathcal Y}$, the connection is the classical
Gauss--Manin connection on the vector bundle over the moduli space of complex
structures on ${\mathcal Y}$, with fibers the mid-dimensional cohomology
$H^{d}({\mathcal Y})$ as mirror symmetry identifies
$H^{k}({\mathcal X}, \wedge ^{k} T^{*}_{\mathcal X})$ with
$H^{k}({\mathcal Y}, \wedge ^{d-k} T^{*}_{\mathcal Y})$.

%s3.2.1 #&#
\subsubsection{}
%%LEAP%%%\label{197.sec3.2.1}
\label{s-solutions}

Solutions to the equation live in a vector space, spanned by K-theory classes
of branes \cite{CV, iritani1, iritani2, KKP}. These are B-type branes on
${\mathcal X}$, objects of ${\mathscr{D}_{\mathcal X}}$, and A-type branes
on ${\mathcal Y}$, objects of $\mathscr{D}_{\mathcal Y}$. A characteristic
feature is that the equation and its solutions mix the A- and B-type structures
on the same manifold.

From the perspective of ${\mathcal X}$, the solutions of the quantum differential
equation come from Gromov--Witten theory. They are obtained by counting
holomorphic maps from a domain curve ${D}$ to ${\mathcal X}$, where
${D}$ is best thought of as an infinite cigar \cite{HIV, mirrorbook} together
with insertions of a class in
$\alpha \in H_{\mathrm{even}}^{*}({\mathcal X})$ at the origin, and
$[{\mathcal F]} \in K({\mathcal X})$ at infinity. The latter is the K-theory
class of a B-type brane ${\mathcal F} \in \mathscr{D}_{\mathcal X}$ which
serves as the boundary condition at the $S^{1}$ boundary at infinity of
${D}$. In the mirror ${\mathcal Y}$, the A- and B-type structures get
exchanged. In the interior of ${D}$, supersymmetry is preserved by
B-type twist, and at the boundary at infinity we place an A-type brane
${\mathcal L} \in \mathscr{D}_{\mathcal Y}$, whose $K$-theory class picks
which solution of the equation we get.

%s3.2.2 #&#
\subsubsection{}
\label{197.sec3.2.2}

One of the key mirror symmetry predictions is that monodromy of the
quantum differential equation gets categorified by the action of derived
autoequivalences of ${\mathscr{D}_{\mathcal X}}$. It is related by mirror symmetry to the monodromy of the
Gauss--Manin connection, computed by Picard--Lefshetz theory, whose categorification
by ${\mathscr{D}_{\mathcal Y}}$ is developed by Seidel \cite{Seidel}.

The flat section ${\mathcal V}$ of the connection in~\eqref{qdif} has a
close cousin. This is Douglas' \cite{Pi1, Pi2, Douglas} $\Pi $-stability
central charge function
${\mathcal Z}^{0}:K({\mathscr D})\rightarrow {\mathbb C}$, whose existence
motivated Bridgeland's formulation of stability conditions \cite{B}. The
$\Pi $-stability central charge ${\mathcal Z}^{0}$ arises from the same
setting as ${\mathcal V}$, except one places trivial insertions at the
origin of ${D}$. This implies that monodromies of ${\mathcal V}$ and
${\mathcal Z}^{0}$ coincide \cite{CV}. In the context of the mirror
${\mathcal Y}$, given any brane
${\mathcal L}\in \mathscr{D}_{\mathcal Y}$, its central charge is simply
${\mathcal Z}^{0}[{\mathcal L}] = \int _{\mathcal L} \Omega $, where
$\Omega $ is the top holomorphic form on ${\mathcal Y}$. The stable objects
are special Lagrangians, on which the phase of $\Omega $ is constant. By
mirror symmetry, monodromy of ${\mathcal Z}^{0}$ is expected to induce
the action of monodromy on ${\mathscr{D}_{\mathcal X}}$. Examples of
braid group actions on the derived categories include works of Khovanov
and Seidel \cite{KS}, Seidel and Thomas \cite{ST}, and others
\cite{TS, B2}.

%s4 #&#
\section{Homological link invariants from B-branes}
\label{197.sec4}

The Knizhnik--Zamolodchikov equation not only has the same flavor as the
quantum differential equation, but for some very special choices of
${\mathcal X}$, they coincide. For the time being, we will take
$^{L}\mathfrak{g}$ to be simply laced, so it coincides with its Langlands
dual ${\mathfrak{g}}$.

%s4.1 #&#
\subsection{The geometry}
%%LEAP%%%\label{197.sec4.1}
\label{ss:coulomb}

The manifold ${\mathcal X}$ may be described as the moduli space of
$G$-monopoles on
%
%e\upshape  4.1 #&#
\begin{equation}
\label{R3}
{\mathbb R}^{3} = {\mathbb R} \times {\mathbb C},
\end{equation}
with prescribed singularities. The monopole group $G$ is related to
$^{L}G$, the Chern--Simons gauge group, by Langlands or electric--magnetic
duality. In Chern--Simons theory, the knots are labeled by representations
of $^{L}G$ and viewed as paths of heavy particles, charged electrically
under $^{L}G$. In the geometric description, the same heavy particles appear
as singular, Dirac-type monopoles of the Langlands dual group $G$. The fact the magnetic description
is what is needed to understand categorification was anticipated by Witten
\cite{WF, WK, WJ, Witten16}.

%s4.1.1 #&#
\subsubsection{}
\label{197.sec4.1.1}

Place a singular $G$ monopole for every finite puncture on
${\mathcal A}\cong {\mathbb R} \times S^{1}$, at the point on
${\mathbb R}$ obtained by forgetting the $S^{1}$. Singular monopole charges
are elements of the cocharacter lattice of $G$, which Langlands duality
identifies with the character lattice of $^{L}G$. Pick the charge of the
monopole to be the highest weight $\mu _{i}$ of the $^{L}G$ representation
$V_{i}$ coloring the puncture. The relative positions of singular monopoles
on ${\mathbb R}^{3}$ are the moduli of the metric on ${\mathcal X}$, so
we will hold them fixed.

The smooth monopole charge is a positive root of $^{L}G$; choose it so
that the total monopole charge is the weight $\nu $ of subspace of representation
$\bigotimes _{i} V_{i}$, where the conformal blocks take values. For our current
purpose, it suffices to assume
%
%e\upshape  4.2 #&#
\begin{equation}
\label{weight}
\nu = \sum _{i} \mu _{i} - \sum _{a=1}^{rk} d_{a} \, ^{L}e_{a},
\end{equation}
is a dominant weight; $^{L}e_{a}$ are the simple positive roots of
$^{L}{\mathfrak{g}}$. Provided $\mu _{i}$ are minuscule co-weights of
$G$ and no pairs of singular monopoles coincide, the monopole moduli space
${\mathcal X}$ is a smooth hyper-Kahler manifold of dimension
\begin{equation*}
\mathrm{dim}_{\mathbb C}({\mathcal X}) =2 \sum _{a} d_{a}.%
\end{equation*}
It is parameterized, in part, by positions of smooth monopoles on
${\mathbb R}^{3}$.

%s4.1.2 #&#
\subsubsection{}
\label{197.sec4.1.2}

A choice of complex structure on ${\mathcal X}$ reflects a split of
${\mathbb R}^{3}$ as ${\mathbb R} \times {\mathbb C}$. The relative positions
of singular monopoles on ${\mathbb C}$ become the complex structure moduli,
and the relative positions of monopoles on ${\mathbb R}$ the Kahler moduli.

This identifies the complexified Kahler moduli space of
${\mathcal X}$ (where the Kahler form gets complexified by a periodic two-form)
with the configuration space of $n$ distinct punctures on
${\mathcal A} = {\mathbb R}\times S^{1}$, modulo overall translations, as in
Figure \ref{f_AR}.

%f2 #&#
\begin{figure}[H]
\begin{center}
\includegraphics{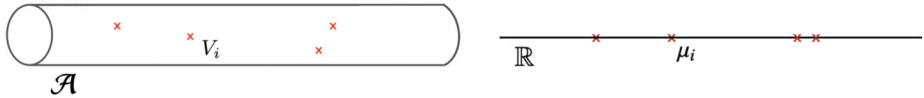}
\end{center}
\caption{Punctures on ${\mathcal A}$ correspond to singular $G$-monopoles
on ${\mathbb R} \in {\mathbb R}\times{\mathbb C}$.}
\label{f_AR}
\end{figure}

%s4.1.3 #&#
\subsubsection{}
\label{197.sec4.1.3}

As a hyper-Kahler manifold, ${\mathcal X}$ has more symmetries than a typical
Calabi--Yau. For its quantum cohomology to be nontrivial, and for the quantum
differential equation to coincide with the KZ equation, we need to work
equivariantly with respect to a torus action that scales its holomorphic
symplectic form
\begin{equation*}
\omega ^{2,0} \rightarrow {\mathfrak{q}} \,\omega ^{2,0}.%
\end{equation*}
For this to be a symmetry, we will place all the singular monopoles at
the origin of ${\mathbb C}$; ${\mathcal X}$ has a
larger torus of symmetries
\begin{equation*}
{T} = {\Lambda} \times {\mathbb C}^{\times}_{\mathfrak{q}},%
\end{equation*}
where ${\Lambda}$ preserves the holomorphic symplectic form, and comes
from the Cartan torus of $G$. The equivariant parameters of the
${\Lambda}$-action correspond to the choice of a flat $^{L}G$ connection
of Chern--Simons theory on ${\mathbb R}^{2}\times S^{1}$.

%s4.1.4 #&#
\subsubsection{}
\label{197.sec4.1.4}

The same manifold ${\mathcal X}$ has appeared in mathematics before, as
a resolution of a transversal slice in the affine Grassmannian
${\mathrm{Gr}}_{G} =G((z))/G[[z]]$ of ${G}$, often denoted by
%
%e\upshape  4.3 #&#
\begin{equation}
\label{Cdef0}
{\mathcal X} = {{{\mathrm{Gr}}}^{{\vec \mu}}}_{\nu}.
\end{equation}
The two are related by thinking of monopole moduli space
${\mathcal X}$ as obtained by a sequence of Hecke modifications of holomorphic
$G$-bundles on ${\mathbb C}$ \cite{KW}.

Manifold ${\mathcal X}$ is also the Coulomb branch of a 3d quiver gauge
theory with ${\mathcal N}=4$ supersymmetry, with quiver based on the Dynkin
diagram of ${\mathfrak{g}}$, see e.g. \cite{GC}. The ranks of
the flavor and gauge symmetry groups are determined from the weights
$\mu $ and $\nu $.

%s4.1.5 #&#
\subsubsection{}
\label{197.sec4.1.5}

The vector ${\vec \mu} =(\mu _{1}, \ldots , \mu _{n})$ in~\eqref{Cdef0} encodes singular monopole charges, and the order in which
they appear on ${\mathbb R}$, and $\nu $ is the total monopole charge.
The ordering of entries of ${\vec \mu}$ is a choice of a chamber in the
Kahler moduli. We will suppress ${\vec \mu}$ for the most part, and denote
all the corresponding distinct symplectic manifolds by
${\mathcal X}$.

%s4.1.6 #&#
\subsubsection{}
\label{197.sec4.1.6}

By a recent theorem of Danilenko \cite{Danilenko}, the Knizhnik--Zamolodchikov
equation corresponding to the Riemann surface
${\mathcal A} = {\mathbb R}\times S^{1}$, with punctures colored by minuscule
representations $V_{i}$ of $^{L}{\mathfrak{g}}$, coincides with the quantum
differential equation of the ${T}$-equivariant Gromov--Witten theory on
${\mathcal X}= {{{\mathrm{Gr}}}^{{\vec \mu}}}_{\nu}$.

This has many deep consequences.

%s4.2 #&#
\subsection{Branes and braiding}
%%LEAP%%%\label{197.sec4.2}
\label{BB}

Since the KZ equation is the quantum-differential equation of ${T}$-equivariant\break
Gromov--Witten theory of ${\mathcal X}$, the space of its solutions gets
identified with $K_{T}({\mathcal X})$, the ${T}$-equivariant $K$-theory
of ${\mathcal X}$.

This is the $K$-group of the category of its B-type branes, the derived
category of ${T}$-equivariant coherent sheaves on ${\mathcal X}$,
\begin{equation*}
{\mathscr D}_{{\mathcal X}} = D^{b}\mathrm{Coh}_{T}({\mathcal X}).
\end{equation*}
This connection between the KZ equation and ${\mathscr D}_{{\mathcal X}} $ is the starting point for our first solution of the categorification problem.

%s4.2.1 #&#
\subsubsection{}
\label{197.sec4.2.1}

A colored braid with $n$ strands in ${\mathcal A} \times [0,1]$ has a geometric
interpretation as a path in the complexified Kahler moduli of
${\mathcal X}$ that avoids singularities, as the order of punctures on
${\mathcal A}$ corresponds to a choice of chamber in the Kahler moduli
of ${\mathcal X}$.

The monodromy of the quantum differential equation along this path acts
on $K_{T}({\mathcal X})$. Since the quantum differential equation coincides
with the KZ equation, by the theorem of \cite{Danilenko},
$K_{T}({\mathcal X})$ becomes a module for
$U_{\mathfrak{q}}(^{L}{\mathfrak{g}})$, corresponding to the weight $\nu$ subspace of representation
$V_{1}\otimes \cdots \otimes V_{n}$.

The fact that derived equivalences of $\mathscr{D}_{\mathcal X}$ categorify
this action is not only an expectation, but also a theorem by Bezrukavnikov
and Okounkov \cite{BO}, whose proof makes use of quantization of
${\mathcal X}$ in characteristic $p$.

%s4.2.2 #&#
\subsubsection{}
\label{197.sec4.2.2}

From physics perspective, the reason derived equivalences of
$\mathscr{D}_{\mathcal X}$ had to categorify the action of monodromy
of the quantum differential equation on $K_{{T}}({\mathcal X})$ is as follows.

Braid group acts, in the $\sigma $-model on the cigar $D$ from Section~\ref{s-solutions}, by letting the moduli of ${\mathcal X}$ vary according
to the braid near the boundary at infinity. The Euclidean time, running
along the cigar, is identified with the time along the braid. This leads
to a Berry phase-type problem studied by Cecotti and Vafa \cite{CV}. It
follows that the $\sigma $-model on the annulus, with moduli that vary
according to the braid, computes the matrix element of the monodromy
$\mathfrak{B}$, picked out by a pair of branes ${\mathcal F}_{0}$ and
${\mathcal F}_{1}$ at its boundaries.

The $\sigma $-model on the same Euclidian annulus, where we take the time
to run around $S^{1}$ instead, computes the index of the supercharge
$Q$ preserved by the two branes. The cohomology of ${Q}$ is computed by
$\mathscr{D}_{\mathcal X}$ as its most basic ingredient, the space of morphisms
\begin{equation*}
{{\mathrm{Hom}}}_{\mathscr{D}_{\mathcal X}}^{*,*}({\mathscr B} {\mathcal F}_{0},
{\mathcal F}_{1})
\end{equation*}
between the branes. This is the space of supersymmetric ground states
of the $\sigma $-model on a strip, obtained by cutting the annulus open.
We took here all the variations of moduli to happen near one boundary,
at the expense of changing a boundary condition from
${\mathcal F}_{0}$ to ${\mathscr B} {\mathcal F}_{0}$. This does not affect
the homology \cite{GMW, A1}, for the very same reason the theory depends on the
homotopy type of the braid only. Per construction, the graded Euler characteristic
of the homology theory, computed by closing the strip back up to the annulus,
is the braiding matrix element,
%
%e\upshape  4.4 #&#
\begin{equation}
\label{EulerX}
\chi ({\mathscr B}\,{\mathcal F}_{0}, {\mathcal F}_{1}) =
({\mathfrak B}\,{\mathcal V}_{0}, {\mathcal V}_{1}),
\end{equation}
between the conformal blocks
${\mathcal V}_{0,1} = {\mathcal V}[{\mathcal F}_{0,1}]$ of the two branes.

Thus, by viewing the same Euclidian annulus in two different ways, we learn
that the braid group action on the derived category
%
%e\upshape  4.5(thesamelabel) \upshape  4.10 #&#
\begin{equation}
\label{BC_}
{\mathscr B}: {\mathscr D}_{{\mathcal X}_{\vec \mu}} \rightarrow
{\mathscr D}_{{\mathcal X}_{\vec \mu '}},
\end{equation}
manifestly categorifies the monodromy matrix
${\mathfrak B}\in U_{\mathfrak{q}}(^{L}\mathfrak{g})$ of the KZ equation.

%s4.3 #&#
\subsection{Link invariants from perverse equivalences}
%%LEAP%%%\label{197.sec4.3}
\label{s_HLI}

The quantum $ U_{\mathfrak{q}}(^{L}\mathfrak{g})$ invariants of knots and
links are matrix elements of the braiding matrix ${\mathfrak B}$, so they
too will be categorified by $\mathscr{D}_{\mathcal X}$, provided we can
identify objects ${\mathcal U} \in {\mathscr{D}_{\mathcal X}}$ which serve
as cups and caps.

Conformal blocks corresponding to cups and caps are defined using fusion
\cite{RCFT}. The geometric analogue of fusion, in terms of
${\mathcal X}$ and its category of branes, was shown in \cite{A1} to be
the existence of certain perverse filtrations on
$\mathscr{D}_{\mathcal X}$, defined by abstractly by Chuang and Rouqiuer
\cite{CR}. The utility of perverse filtrations for understanding the action
of braiding on $\mathscr{D}_{\mathcal X}$ parallels the utility of fusion
in describing the action of braiding in conformal field theory. In particular,
it leads to identification of the cup and cap branes ${\mathcal U}$ we
need, and a simple proof that
${{\mathrm{Hom}}}_{\mathscr{D}_{\mathcal X}}^{*,*}({\mathscr B} {\mathcal U},
{\mathcal U})$ are homological invariants of links \cite{A1}.

%s4.3.1 #&#
\subsubsection{}
\label{197.sec4.3.1}

As we bring a pair of punctures at $a_{i}$ and $a_{j}$ on
${\mathcal A}$ together, we get a new natural basis of solutions to the
KZ equation, called the fusion basis, whose virtue is that it diagonalizes
braiding. The possible eigenvectors are labeled by the representations
%
%e\upshape  4.6 #&#
\begin{equation}
\label{tensor}
{V}_{i} \otimes {V}_{j} =\bigotimes _{m=0}^{m_{\max}} {V}_{k_{m}},
\end{equation}
that occur in the tensor product of representations $V_{i}$ and
$V_{j}$ labeling the punctures. Because $V_{i}$ and $V_{j}$ are minuscule
representations, the nonzero multiplicities on
the right-hand side are all equal to $1$. The cap arises as a special case,
obtained by starting with a pair of conjugate representations
$V_{i}$ and $V_{i}^{\star}$, and picking the trivial representation in their
tensor product.

%s4.3.2 #&#
\subsubsection{}
\label{197.sec4.3.2}

From perspective of ${\mathcal X}$, a pair
of singular monopoles of charges $\mu_i$ and $\mu_j$ are coming together on ${\mathbb R}$, as in Figure \ref{f_AR}, and we approach a wall in Kahler moduli
at which ${\mathcal X}$ develops a singularity. At the singularity, a collection
of cycles vanishes.
This is due to monopole bubbling phenomena described by Kapustin and Witten
in \cite{KW}.

The types of monopole bubbling that can occur
are labeled by representations $V_{k_m}$ that occur in the tensor product
$V_{i} \otimes V_{j}$. The moduli space of monopoles whose positions we
need to tune for the bubbling of type $V_{k_m}$ to occur is
${\mathrm{Gr}}^{(\mu _ {i}, \mu _{j})}_{\mu _{k_m}}= T^{*}F_{k_m}$, where
$\mu _{k_m}$ is the highest weight of $V_{k_m}$. This space is transverse
to the locus where exactly $\mu _{i}+\mu _{j} -\mu _{k_m}$ monopoles have
bubbled off \cite{A1}. It has a vanishing cycle $F_{k_m}$, corresponding
to the representation $V_{k_m}$, as its zero section. (Viewing
${\mathcal X}$ as the Coulomb branch, monopole bubbling is related to partial
Higgsing phenomena.)

%s4.3.3 #&#
\subsubsection{}
\label{197.sec4.3.3}

Conformal blocks which diagonalize the action of braiding do not in general
come from actual objects of the derived category
$\mathscr{D}_{\mathcal X}$. As is well known from Picard--Lefshetz theory,
eigensheaves of braiding, branes on which the braiding acts only by degree
shifts
${\mathscr B} {\mathcal E} = {\mathcal E}[D_{\mathcal E}]\{C_{\mathcal E}\}$, are very rare.

What one gets instead \cite{A1} is a filtration
%
%e\upshape  4.7 #&#
\begin{equation}
\label{filtration0a}
{{\mathscr D}}_{k_{0}} \subset {{\mathscr D}}_{k_{1}} \subset \dots
\subset {{\mathscr D}}_{k_{\max}} = {{\mathscr D}}_{\mathcal X},
\end{equation}
by the order of vanishing of the $\Pi $-stability central charge
${\mathcal Z}^{0}: K({\mathcal X}) \rightarrow {\mathbb C}$. More precisely,
one gets a pair of such filtrations, one on each side of the wall. Crossing
the wall preserves the filtration, but has the effect of mixing up branes
at a given order in the filtration, with those at lower orders, whose central
charge vanishes faster. Because ${\mathcal X}$ is hyper-Kahler, the
$\Pi $-stability central charge is given in terms of classical geometry
(by Eq. (4.7) of \cite{A1}).

The existence of the filtration with the stated properties follows from
the existence of the {equivariant} central charge function
${\mathcal Z,}$
%
%e\upshape  4.8 #&#
\begin{equation}
\label{CF}
{\mathcal Z}: K_{T}({\mathcal X}) \rightarrow {\mathbb C},
\end{equation}
and the fact the action of braiding on $K_{T}({\mathcal X})$ lifts to
the action on $\mathscr{D}_{\mathcal X}$, by the theorem of \cite{BO}.
The equivariant central charge ${\mathcal Z}$ is computed by the equivariant
Gromov--Witten theory on ${\mathcal X}$ in a manner analogous to
${\mathcal V}$, starting with the $\sigma $-model on the cigar
${D}$ except with no insertion at its tip. It reduces to the $\Pi $-stability
central charge ${\mathcal Z}^{0}$ by turning the equivariant parameters
off.

%s4.3.4 #&#
\subsubsection{}
\label{197.sec4.3.4}

While ${\mathscr B}$ has few eigensheaves in
${\mathscr{D}_{\mathcal X}}$, it acts by degree shifts
%
%e\upshape  4.9 #&#
\begin{equation}
\label{BA}
{\mathscr B}: {{\mathscr D}}_{k_{m}}/{{\mathscr D}}_{k_{m-1}}
\rightarrow {{\mathscr D}}_{k_{m}}/{{\mathscr D}}_{k_{m-1}}[ D_{k_{m}}]
\{C_{k_{m}}\},
\end{equation}
on the quotient subcategories. The degree shifts may be read off from the
eigenvectors of the action of braiding on the equivariant central charge
function ${\mathcal Z}$. As the punctures at $a_{i}$ and $a_{j}$ come together,
the eigenvector corresponding to the representation $V_{k_m}$ in \eqref{tensor}, vanishes as \cite{A1}
\begin{equation*}
{\mathcal Z}_{k_{m}} = (a_{i} - a_{j})^{D_{k_{m}} + C_{k_{m}}/\kappa}
\times \textup{finite}.
\end{equation*}
It follows that braiding $a_{i}$ and $a_{j}$ counterclockwise acts
by
\begin{equation*}
{\mathcal Z}_{k_{m}}\rightarrow (-1)^{D_{k_{m}}}
{\mathfrak{q}}^{{\frac{1}{2}} C_{k_{m}}} {\mathcal Z}_{k_{m}}.%
\end{equation*}
The cohomological degree shift
$D_{k_{m}} = \mathrm{dim}_{\mathbb C}\, F_{k_{m}}$ is by the dimension
of the vanishing cycle. The equivariant degree shift $C_{k_{m}}$ is essentially
the one familiar from the action of braiding on conformal blocks of
$\widehat{^L\mathfrak{g}}$ in the fusion basis \cite{A1}.

%s4.3.5 #&#
\subsubsection{}
\label{197.sec4.3.5}

The derived equivalences of this type are the perverse equivalences of
Chuang and Rouquier \cite{CR0, CR}. They envisioned them as a way to describe
derived equivalences which come from variations of Bridgeland stability
conditions, but with few examples from geometry.

Traditionally, braid group actions on derived categories of coherent sheaves,
or B-branes, are fairly difficult to describe, see for example
\cite{CK1, CK2}. Braid group actions on the categories of A-branes are
much easier to understand, via Picard--Lefshetz theory and its categorical
uplifts \cite{Seidel}, see e.g. \cite{KS, TS}. The theory
of variations of stability conditions, by Douglas and Bridgeland, was invented
to bridge the two \cite{Pi1, Douglas}.

%s4.3.6 #&#
\subsubsection{}
%%LEAP%%%\label{197.sec4.3.6}
\label{s:Gras}

As a by-product, we learn that conformal blocks describing collections
cups or caps colored by minuscule representations, come from branes in $\mathscr{D}_{\mathcal X}$ which
have a simple geometric meaning \cite{A1}.

Take
${\mathcal X}= \mathrm{Gr}^{{(\mu _{1}, \mu _{1}^{*}, \ldots , {\mu _{d}, \mu _{d}^{*}})}}_{0}$
corresponding to ${\mathcal A}$ with $n=2d$ punctures, colored by pairs of complex conjugate, minuscule  representations $V_{i}$ and
$V_{i}^{*}$. We get a vanishing cycle $U$ in ${\mathcal X}$ which is
a product of $d$ minuscule Grassmannians,
\begin{equation*}
U = G/P_{1} \times \dots \times G/P_{d},
\end{equation*}
where $P_{i}$ is the maximal parabolic subgroup of $G$ associated to representation
$V_{i}$. This vanishing cycle embeds in ${\mathcal X}$ as a compact holomorphic
Lagrangian, so in the neighborhood of $U$, we can model
${\mathcal X}$ as $T^{*}U$. The structure sheaf
\begin{equation*}
{\mathcal U} = {\mathcal O}_{U} \in \mathscr{D}_{\mathcal X}
\end{equation*}
of $U$ is the brane we are after. The Grassmannian
$G/P_{i}$ is the cycle that vanishes when a single pair of singular monopoles
of charges $\mu _{i}$ and $\mu _{i}^{*}$ come together, as
${\mathrm{Gr}}^{(\mu _ {i}, \mu _{i}^{*})}_{0} = T^{*}G/P_{i}$.

The brane ${\mathcal U}$ lives at the very bottom of a $d$-fold filtration
which $\mathscr{D}_{\mathcal X}$ develops at the intersection of $d$ walls in the
Kahler moduli of ${\mathcal X}$ corresponding to bringing punctures
together pairwise. It follows ${\mathcal U}$ is the eigensheaf of braiding
each pair of matched endpoints. It is extremely special, for the same reason
the trivial representation is special.

%s4.3.7 #&#
\subsubsection{}
\label{197.sec4.3.7}

Just as fusion provides the right language to understand the action of
braiding in conformal field theory, the perverse filtrations provide the right
language to describe the action of braiding on derived categories.
Using perverse filtrations and the very special properties of the vanishing
cycle branes ${\mathcal U} \in \mathscr{D}_{\mathcal X}$, one gets the
following theorem \cite{A1}:
%
%t1 #&#
\begin{theorem}%
\label{one}
For any simply laced Lie algebra $^{L}{\mathfrak{g}}$, the homology groups
\begin{equation*}
\mathrm{Hom}^{*,*}_{\mathscr{D}_{\mathcal X}}({\mathscr B} {\mathcal U}, {\mathcal U}),
\end{equation*}
categorify $U_{\mathfrak{q}}(^{L}{\mathfrak{g}})$ quantum link invariants,
and are themselves link invariants.
\end{theorem}

%s4.3.8 #&#
\subsubsection{}
\label{197.sec4.3.8}

As an illustration, proving that (the equivalent of) the pitchfork move in
the figure below holds in ${\mathscr{D}_{\mathcal X}}$
%
%f3 #&#
\begin{figure}[H]
\begin{center}
\includegraphics{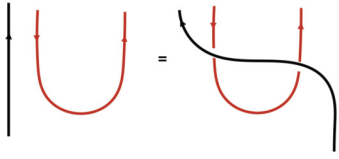}
\end{center}
\caption{A move equivalent to the pitchfork move.}
\label{f_rel1}
\end{figure}
\noindent
requires showing that we have a derived equivalence
%
%e\upshape  4.10(thesamelabel) \upshape  4.5 #&#
\begin{equation}
\label{BC}
{\mathscr B}\circ {\mathscr C}_{i} \cong {\mathscr C}_{i}'',
\end{equation}
where ${\mathscr C}_{i}$ and ${\mathscr C}_{i}''$ are cup functors on the
right and the left in Figure~\ref{f_rel1}, respectively. They increase
the number of strands by two and map
\begin{equation*}
{\mathscr C}_{i} : {\mathscr D}_{{\mathcal X}_{n-2}}  \rightarrow {
\mathscr D}_{{\mathcal X}_{n}} \quad \textup{and}\quad {\mathscr C}''_{i} :
{\mathscr D}_{{\mathcal X}_{n-2}} \rightarrow {\mathscr D}_{{\mathcal X}''_{n}},%
\end{equation*}
where the subscript serves to indicate the number of strands. The functor
${\mathscr B}$ is the equivalence of categories from the theorem of
\cite{BO}
\begin{equation*}
{\mathscr B}: {\mathscr D}_{{\mathcal X}_{n}} \rightarrow {\mathscr D}_{{\mathcal X}''_{n}},%
\end{equation*}
corresponding to braiding ${V}_{k}(a_{k})$ with
${V_{i}}(a_{i}) \otimes {{V}_{i}^{*}}(a_{j})$ where $V_{i}$ and
$V_{i}^{*}$ color the red and $V_{k}$ the black strand in Figure~\ref{f_rel1}.

To prove the identity~\eqref{BC} note that
%
%e\upshape  4.11 #&#
\begin{equation}
\label{gra}
{\mathscr C}_{i} {\mathscr D}_{{\mathcal X}_{n-2}} \subset
{\mathscr D}_{{\mathcal X}_{n}} \quad \textup{and}\quad  {\mathscr C}''_{i}
{\mathscr D}_{{\mathcal X}_{n-2}} \subset {\mathscr D}_{{\mathcal X}_{n}''},
\end{equation}
are the subcategories which are the bottom-most part of the double filtrations of
$ {\mathscr D}_{{\mathcal X}_{n}}$ and
$ {\mathscr D}_{{\mathcal X}_{n}''}$, corresponding to the intersection of
walls at which the three punctures come together. By the definition
of perverse filtrations, the functor ${\mathscr B}$ acts at a bottom part
of a double filtration at most by degree shifts. The degree shifts are
trivial too, since if they were not, the relation we are after would not
hold even in conformal field theory, and we know it does. To complete the
proof, one recalls that a perverse equivalence that acts by degree shifts
that are trivial is an equivalence of categories \cite{CR}.

Proofs of invariance under the Reidermeister 0 and the framed Reidermeister
I moves are similar. The invariance under Reidermeister II and III moves
follows from the theorem of \cite{BO}. One should compare this to proofs
of the same relations in \cite{CK1,CK2}, which are more technical and less
general.

%s4.3.9 #&#
\subsubsection{}
\label{197.sec4.3.9}

An elementary consequence is a geometric explanation of mirror symmetry
which relates the $U_{\mathfrak{q}}(^{L}\mathfrak{g})$ invariants of a
link $K$ and its mirror reflection $K^{*}$.

It is a consequence of a basic property of
${\mathscr{D}_{\mathcal X}}$, Serre duality. Serre duality implies the
isomorphism of homology groups on ${\mathcal X}$ which is a $2d$ complex-dimensional
Calabi--Yau manifold,
%
%e\upshape  4.12 #&#
\begin{equation}
\label{mirr3}
{\mathrm{Hom}}_{{\mathscr D}_{\mathcal X}}\bigl({\mathscr B} \,{\mathcal U},
{\mathcal U}[M]\{J_{0}, {\vec J}\}\bigr)
={\mathrm{Hom}}_{{\mathscr D}_{\mathcal X}}
\bigl({\mathscr B}{\mathcal U}, {\mathcal U}[2d-M]\{-d-J_{0}, -{\vec J}\}\bigr).
\end{equation}
The equivariant degree shift comes from the fact the unique holomorphic section of the canonical bundle has weight $d$ under the
${\mathbb C}^{\times}_{\mathfrak{q}}\subset {\mathrm T}$ action.  Mirror symmetry follows by taking Euler
characteristic of both sides \cite{A1}.

%s4.4 #&#
\subsection{Algebra from B-branes}
%%LEAP%%%\label{197.sec4.4}
\label{s-BK}

Bezrukavnikov and Kaledin, using quantization in characteristic $p$, constructed
a tilting vector bundle ${\mathcal T}$, on any ${\mathcal X}$ which is
a symplectic resolution \cite{kaledin, kaledinR, BK1, BK2}. Its endomorphism
algebra
\begin{equation*}
{\mathscr A} = \mathrm{Hom}^{*}_{\mathscr{D}_{\mathcal X}}({\mathcal T}, {\mathcal T})
\end{equation*}
is an ordinary associative algebra, graded only by equivariant degrees.
The derived category $\mathscr{D}_{\mathscr A}$ of its modules is equivalent
to ${\mathscr{D}_{\mathcal X}}$,
\begin{equation*}
\mathscr{D}_{\mathcal X}\cong \mathscr{D}_{\mathscr A},%
\end{equation*}
essentially per definition.

Webster recently computed the algebra ${\mathscr A}$ for our
${\mathcal X}$ \cite{W2}, and showed that it coincides with a cylindrical
version of the KLRW algebra from \cite{webster}. Working with the cylindrical
KLRW algebra, as opposed to the ordinary one, leads to invariants of links
in ${\mathbb R}^{2}\times S^{1}$ and not just in ${\mathbb R}^{3}$. The
KLRW algebra generalizes the algebras of Khovanov and Lauda \cite{KL_} and
Rouquier \cite{R}. The cylindrical version of the KLR algebra corresponds
to ${\mathcal X}$ which is a Coulomb branch of a pure 3D gauge
theory.

%s4.4.1 #&#
\subsubsection{}
\label{197.sec4.4.1}

The description of link homologies via
$\mathscr{D}_{\mathcal X}= D\!\operatorname{Coh}_{T}({\mathcal X})$ provides a geometric
meaning of homological $U_{\mathfrak{q}}(^{L}{\mathfrak{g}})$ link invariants.
Even so, without further input, the description of link homologies either
in terms of $\mathscr{D}_{\mathcal X}$ or $\mathscr{D}_{\mathscr A}$ is
purely formal. With the help of (equivariant) homological mirror symmetry,
we will give a description of link homology groups which is
explicit and explicitly computable; in this sense, link homology groups
come to life in the mirror.

%s5 #&#
\section{Mirror symmetry for monopole moduli space}
\label{197.sec5}

In the very best instances, homological mirror symmetry relating
${\mathscr{D}_{\mathcal Y}}$ and ${\mathscr{D}_{\mathcal X}}$ can be made
manifest, by showing that each is equivalent to
$\mathscr{D}_{\mathscr A}$, the derived category of modules of the same
associative algebra ${\mathscr A,}$
%
%e\upshape  5.1 #&#
\begin{equation}
\label{XAYa}
{\mathscr{D}_{\mathcal X}} \cong {\mathscr{D}_{\mathscr A}} \cong {\mathscr{D}_{\mathcal Y}}.
\end{equation}
The algebra
\begin{equation*}
{\mathscr A}=\mathrm{Hom}^{*}_{\mathscr D}({\mathcal T}, {\mathcal T})%
\end{equation*}
is the endomorphism algebra of a set of branes
${\mathcal T} = \bigoplus _{{\mathcal C}} {\mathcal T}_{\mathcal C}$, which
generate ${\mathscr{D}_{\mathcal X}}$ and
${\mathscr{D}_{\mathcal Y}}$. For economy, we will be denoting branes related by mirror symmetry by the same letter.

An elementary example \cite{Auroux} is mirror symmetry relating a pair
of infinite cylinders, ${\mathcal X}= {\mathbb C}^{\times}$ and
${\mathcal Y}= {\mathbb R}\times S^{1}$, whose torus fibers are dual $S^{1}$'s. Both ${\mathscr{D}_{\mathcal X}}$, the derived category
of coherent sheaves on ${\mathcal X}$, and
$\mathscr{D}_{\mathcal Y}$, based on the wrapped Fukaya category, are generated
by a single object ${\mathcal T}$, a flat line bundle on
${\mathcal X}$ and a real-line Lagrangian on ${\mathcal Y}$. Their algebras
of open strings are the same, equal to the algebra
${\mathscr A} = {\mathbb C}[x^{\pm 1}]$ of holomorphic functions on the
cylinder.

%f4 #&#
\begin{figure}[H]
\vspace*{-6pt}
\begin{center}
\includegraphics{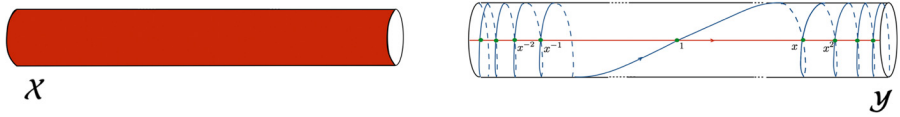}
\end{center}
\caption{A simple example of manifest mirror symmetry.}
\label{f_cyl}
\vspace*{-15pt}
\end{figure}

\enlargethispage{3pt}
%s5.1 #&#
\subsection{The algebra for homological mirror symmetry}
%%LEAP%%%\label{197.sec5.1}
\label{s-mirrorm}

In our setting, the generator ${\mathcal T}$ of
$\mathscr{D}_{\mathcal X}$ is the tilting generator of Bezrukavnikov and
Kaledin from Section~\ref{s-BK}. Webster's proof of the equivalence of
categorification of $U_{\mathfrak{q}}(^{L}{\mathfrak{g}})$ link invariants
and B-type branes on ${\mathcal X}$ and via the cKLRW algebra
${\mathscr A}$ should be understood as the first of the two equivalences
in~\eqref{XAYa}.

\vspace*{-3pt}
%s5.1.1 #&#
\subsubsection{}
\label{197.sec5.1.1}

The mirror ${\mathcal Y}$ of ${\mathcal X}$ is the moduli space of
$G$ monopoles, of the same charges as ${\mathcal X}$ except on
${\mathbb R}^{2}\times S^{1}$ instead of on ${\mathbb R}^{3}$, with only complex
and no Kahler moduli turned on, and equipped with a potential
\cite{A2}. Without the potential, the mirror to ${\mathcal Y}$ would be
another moduli space of $G$ monopoles on
${\mathbb R}^{2}\times S^{1}$.

The theory based on $\mathscr{D}_{\mathcal Y}$, the derived Fukaya--Seidel
category of ${\mathcal Y}$, is in the same spirit as the work of Seidel and Smith \cite{SSm}. They pioneered geometric approaches to link homology,
but produced a only singly graded theory, known as symplectic Khovanov homology.
The computation of $\mathscr{D}_{\mathcal Y}$, which makes mirror symmetry in~\eqref{XAYa} manifest, is given in the joint work with Danilenko, Li, and
Zhou \cite{ADZ}.

\vspace*{-3pt}
%s5.2 #&#
\subsection{The core of the monopole moduli space}
\label{197.sec5.2}

Working equivariantly with respect to a
${\mathbb C}_{\mathfrak{q}}^{\times}$-symmetry which scales the holomorphic
symplectic form of ${\mathcal X}$, all the information about its geometry
should be encoded in a core locus preserved by such actions.

The core $X$ is a singular holomorphic Lagrangian in ${\mathcal X}$ which
is the union of supports of all stable envelopes \cite{MO, ese}. Equivalently,
$X$ is the union of all attracting sets of $\Lambda $-torus actions on
${\mathcal X}$, where we let $\Lambda $ vary over all chambers. If we view
${\mathcal X}$ as the monopole moduli space, we can put this more simply:
$X$ is the locus where all the monopoles, singular or not, are at the origin
of ${\mathbb C}$ in ${\mathbb R} \times {\mathbb C}$. Viewing it as a Coulomb
branch, $X$ is the locus at which the complex scalar fields in vector multiplets
vanish.

We will define the equivariant mirror $Y$ of ${\mathcal X}$ to be the ordinary
mirror of its core, so we have
\begin{center}
\begin{tikzcd}[sep= large]
{\mathcal X}\ar[r,<->,"\textup{mirror}"]
\arrow[<->]{rd}[description]{\textup{equiv. mirror}}&{\mathcal Y}
\ar[d, <->,""]
\\
X\ar[r,<->,"\textup{mirror}"]\ar[u,<->,""] & Y
\end{tikzcd}
\end{center}%
\label{core2}
\vskip 0.15cm
Working equivariantly with respect to the ${T}$-action on
${\mathcal X}$, the equivariant mirror gets a potential $W$, making the
theory on $Y$ into a Landau--Ginsburg model. While $X$ embeds into
${\mathcal X}$ as a holomorphic Lagrangian of dimension $d$,
${\mathcal Y}$ fibers over $Y$ with holomorphic Lagrangian
$({\mathbb C}^{\times})^{d}$ fibers.

%s5.2.1 #&#
\subsubsection{}
%%LEAP%%%\label{197.sec5.2.1}
\label{AIS}

A model example is ${\mathcal X}$ which is the resolution of an
$A_{n-1}$ hypersurface singularity, $uv=z^{n}$; ${\mathcal X}$ is the moduli
space of a single smooth $G= \mathrm{SU}(2)/{\mathbb Z}_{2}$ monopole, in the presence
of $n$ singular ones. The core $X$ is a collection of $n-1$
${\mathbb P}^{1}$'s with a pair of infinite discs attached, as in Figure \ref{f_core}.
%
%f5 #&#
\begin{figure}[h!]
\begin{center}
\includegraphics{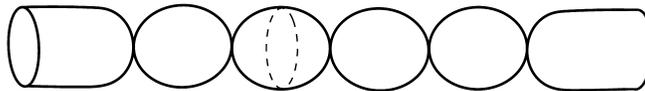}
\end{center}
\caption{Core $X$ of a resolution of the $A_{n-1}$ singularity.}
\label{f_core}
\end{figure}

The ordinary mirror ${\mathcal Y}$ of ${\mathcal X}$ is the complex structure
deformation of the ``multiplicative'' $A_{n-1}$ surface singularity, with
a potential which we will not need. ${\mathcal Y}$ is a
${\mathbb C}^{\times}$ fibration over $Y$ which is itself an infinite cylinder,
%
%f6 #&#
\begin{figure}[!hbtp]
\begin{center}
\includegraphics{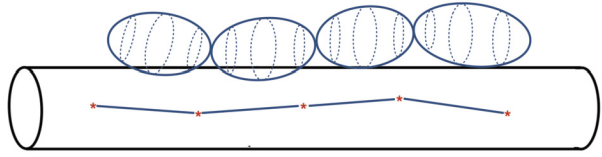}
\end{center}
\caption{Lagrangian spheres in ${\mathcal Y}$ mirror the vanishing
${{\mathbb P}^{1}}$'s in ${\mathcal X}$.}
\label{f_Ycore}
\end{figure}
a copy of ${\mathbb C}^{\times}$ with $n$ points deleted. At the marked
points, the ${\mathbb C}^{\times}$ fibers degenerate. There are
$n-1$ Lagrangian spheres in ${\mathcal Y}$, which are mirror to
$n-1$ ${\mathbb P}^{1}$'s in ${\mathcal X}$. They project to Lagrangians
in $Y$ which begin and end at the punctures.

%s5.2.2 #&#
\subsubsection{}
\label{197.sec5.2.2}

The model example corresponds to $^{L}G = \mathrm{SU}(2)$ Chern--Simons theory on
${\mathbb R}^{2}\times S^{1}$, and $\widehat{{\mathfrak {su}}_{2}}$ conformal
blocks on ${\mathcal A} ={\mathbb R}\times S^{1}$. The $n$ punctures on
${\mathcal A}$ are colored by the fundamental, two-dimensional representation
$V_{1/2}$ of ${\mathfrak{su}}_{2}$, and we take the subspace of weight
one level below the highest. Note that $Y$ coincides with the Riemann
surface ${\mathcal A}$ where the conformal blocks live. This is not an
accident.

In the model example, both $X$ and $Y$ are $S^{1}$ fibrations over
${\mathbb R}$ with $n$ marked points. At the marked points, the
$S^{1}$ fibers of $X$ degenerate. In $Y$, this is mirrored by fibers that decompactify,
due to points which are deleted.

%s5.2.3 #&#
\subsubsection{}
\label{197.sec5.2.3}

More generally, for ${\mathcal X} = {\mathrm{Gr}^{\vec \mu}}_{\nu}$ we have
$d_{a}$ smooth $G$-monopoles colored by simple roots $^{L}e_{a}$ and otherwise
identical. It follows that the common base of SYZ fibrations of $X$ and
$Y$ is the configuration space of the smooth monopoles on the real line
${\mathbb R}$ with $n$ marked points. The marked points are labeled by the weights
$\mu _{i}$ of $^{L}{\mathfrak{g}}$, which are the singular monopole charges.

An explicit description of $Y$, as well as its category of A-branes
$\mathscr{D}_{Y}$, is given \cite{ADZ}. Here we will only describe some
of its features. In an open set, $Y$ coincides with
\begin{equation*}
Y_{0} =\bigotimes _{a =1}^{rk} \mathrm{Sym}^{d_{a}}{\mathcal A},%
\end{equation*}
the configuration space of $d=\sum _{a=1}^{rk} d_{a}$ points on the punctured
Riemann surface ${\mathcal A}$, ``colored'' by simple roots
$^{L}e_{a}$ of $^{L}{\mathfrak{g}}$, but otherwise identical. The open set
is the complement of the divisor of zeros and of poles of function
$f^{0}$ in~\eqref{dressed}.

The top holomorphic form on $Y$ is
%
%e\upshape  5.2 #&#
\begin{equation}
\label{Omega}
\Omega = \bigwedge _{a=1}^{rk} \bigwedge _{\alpha =1}^{d_{a}} {\frac{dy_{\alpha , a}}{y_{\alpha , a}}},
\end{equation}
where $ y_{\alpha , a}$ are coordinates on $d$ copies of
${\mathcal A}$, viewed as the complex plane with $0$ and $\infty $ deleted.
While $\Omega $ itself is not globally well defined, so $K_{Y}$ is not
trivial, $\Omega ^{\otimes 2}$ is well defined and
%
%e\upshape  5.3 #&#
\begin{equation}
\label{2c}
2c_{1}(K_{Y}) =0.
\end{equation}
This allows $\mathscr{D}_{Y}$ to have a ${\mathbb Z}$-valued cohomological
grading. The symplectic form on $Y$ is inherited from the symplectic form
on ${\mathcal Y}$, by restricting it to the vanishing $(S^{1})^{d}$ in
each of its $({\mathbb C}^{\times})^{d}$ fibers over $Y$ \cite{ADZ}. The
precise choice of symplectic structure is the one compatible with mirror
symmetry which we used to define $Y$, as the equivariant mirror of
${\mathcal X} = {\mathrm{Gr}^{\vec \mu}}_{\nu}$ and the ordinary mirror of its core.

Including the equivariant ${T}$-equivariant action on
${\mathcal X}$ and $X$ corresponds to adding to the $\sigma $-model on
$Y$ a potential
%
%e\upshape  5.4 #&#
\begin{equation}
\label{sup}
W = \lambda _{0} W^{0} + \sum _{a=1}^{rk} \,\lambda _{a}\,W^{a},
\end{equation}
which is a multivalued holomorphic function on $Y$; $\lambda _{a}$ are
the equivariant parameters of the ${\Lambda}$-action on
${\mathcal X}$, and
\begin{equation*}
{\mathfrak{q}} = e^{2\pi i \lambda _{0}}.
\end{equation*}
The potentials $W^{0}$ and $W^{a}$ are given by
\begin{equation*}
W^{0} = \ln f^{0}, \quad W^{a} =\ln \prod _{\alpha =1}^{d_{a}}y_{a,
\alpha},
\end{equation*}
where
%
%e\upshape  5.5 #&#
\begin{equation}
\label{dressed}
f^{0} (y)= \prod _{a=1}^{rk} \prod _{\alpha =1}^{d_{a}} {
\frac{\prod _{i}(1- a_{i}/y_{\alpha , a})^{\langle ^{L}e_{a},
\mu _{i} \rangle}}{\prod _{(b,\beta ) \neq (a, \alpha )}
(1- y_{\beta ,b}/y_{\alpha , a})^{\langle ^{L}e_{a}, ^{L}e_{b}\rangle /2}}}.
\end{equation}
The superpotential $W$ breaks the conformal invariance of the
$\sigma $-model to $Y$ if $\lambda _{0} \neq 0$, since only a quasihomogenous
superpotential is compatible with it. This is mirror to breaking of
conformal invariance on ${\mathcal X}$ by the
${\mathbb C}^{\times}_{\mathfrak{q}}$-action for
${\mathfrak{q}}\neq 1$.

Since $W^{0}$ and $W^{a}$ are multivalued, $Y$ is equipped with a collection
of closed one-forms with integer periods
\begin{equation*}
c^{0} = d\!W^{0}/2\pi i,\quad c^{a}=d\!W^{a}/2\pi i \in H^{1}(Y,
\mathbb Z),
\end{equation*}
which introduce additional gradings in the category of A-branes, as in
\cite{SS}.

%s5.2.4 #&#
\subsubsection{}
\label{197.sec5.2.4}

From the mirror perspective, the conformal blocks of
$\widehat{^L\mathfrak{g}}$ come from the B-twisted Landau--Ginsburg model
$(Y,W)$ on ${D}$ which is an infinitely long cigar, with A-type boundary
condition at infinity corresponding to a Lagrangian $L\in Y$. The partition
function of the theory has the following form:

%e\upshape  5.6 #&#
\begin{equation}
\label{VLG}
{\mathcal V}_{\alpha}[L] = \int _{L} \,{\Phi}_{\alpha }\, \Omega \, e^{-W}
,%
\end{equation}
where $\Phi _{\alpha}$ are chiral ring operators, inserted at the tip of
the cigar \cite{CV, HIV, mirrorbook}. By placing the trivial insertion
at the origin instead, we get the equivariant central charge function
${\mathcal Z}[L] = \int _{L} \Omega \,e^{-W }$; by further turning the
equivariant parameters off, the potential $W$ vanishes and the equivariant
central charge becomes the ordinary brane central charge
${\mathcal Z}^{0}[L] = \int _{L} \Omega $.

We have (re)discovered, from mirror symmetry, an integral representation
of the conformal blocks of $\widehat{^L\mathfrak{g}}$. This ``free field
representation'' of conformal blocks, remarkable for its simplicity
\cite{EFK}, goes back to the 1980s work of Kohno and Feigin and Frenkel
\cite{Koh, FF}, and of Schechtman and Varchenko \cite{SV1,SV2}.

%s5.2.5 #&#
\subsubsection{}
\label{197.sec5.2.5}

There is a reconstruction theory, due to Givental \cite{G2} and Teleman
\cite{T}, which says that, starting with the solution of the quantum differential
equation or its mirror counterpart, one gets to reconstruct all genus topological
string amplitudes for any semi-simple 2D field theory. The semi-simplicity condition
is satisfied in our case, as $W$ has isolated critical points. It follows
the B-twisted Landau--Ginsburg model on $(Y, W)$ and A-twisted
${T}$-equivariant sigma model on ${\mathcal X}$ are equivalent to all
genus \cite{A2}. Thus, equivariant mirror symmetry holds as an equivalence
of topological string amplitudes.

%s5.3 #&#
\subsection{Equivariant Fukaya--Seidel category}
\label{197.sec5.3}

For every A-brane $L$ at the boundary at infinity of the cigar ${D}$ ,
we get a solution of the KZ equation. The brane is an object of
\begin{equation*}
\mathscr{D}_{Y}=D\bigl({\mathcal FS}(Y, W)\bigr),%
\end{equation*}
the derived Fukaya--Seidel category of $Y$ with potential $W$. The category
should be thought of as a category of equivariant A-branes, due to the
fact $W$ in~\eqref{sup} is multivalued. Another novel aspect of
$\mathscr{D}_{Y}$ is that it provides an example of Fukaya--Seidel category
with coefficients in perverse schobers. This structure, inherited from
equivariant mirror symmetry, was discovered in \cite{ADZ}.

%s5.3.1 #&#
\subsubsection{}
\label{197.sec5.3.1}

Objects of $\mathscr{D}_{Y}$ are Lagrangians in $Y$, equipped with some
extra data. A Lagrangian in $Y$ is a product of $d$ one-dimensional curves
on ${\mathcal A}$ which are colored by simple roots and may be immersed;
or a simplex obtained from an embedded curve, as a configuration space of  $d$ partially ordered colored points. The theory also includes
more abstract branes, which are iterated mapping cones over morphisms
between Lagrangians.

%s5.3.2 #&#
\subsubsection{}
\label{197.sec5.3.2}

The extra data includes a grading by Maslov and equivariant degrees. The
equivariant grading of a brane in $\mathscr{D}_{Y}$ is defined by choosing
a lift of the phase of $e^{-W}$ to a real-valued function on the Lagrangian
$L$. The equivariant degree shift operation,
\begin{equation*}
L \rightarrow L\{\vec d\},
\end{equation*}
with $\vec{d}\in {\mathbb Z}^{{rk}+1}$, corresponds to changing the lift
of $W$ on $L$, now viewed as a graded Lagrangian,
$ W|_{{L}\{\vec d\}} = W|_{L} + 2\pi i {\vec \lambda}\cdot {\vec d}
$. This is analogous to how a choice of a lift of the phase of
$ \Omega ^{\otimes 2}$ defines the Maslov, or cohomological, grading of
a Lagrangian. This restricts the Lagrangians that give rise to objects
of $\mathscr{D}_{Y}$ to those for which such lifts can be defined.

More generally, branes in $\mathscr{D}_{Y}$ are graded Lagrangians
$L$ equipped with an extra structure of a local system $\Lambda $ of modules
of a certain algebra ${\mathcal B}$ we will describe shortly. For the time
being, only branes for which $\Lambda $ is trivial will play a role for
us.

%s5.3.3 #&#
\subsubsection{}
\label{197.sec5.3.3}

The space of morphisms between a pair of Lagrangian branes in a derived
Fukaya category
\begin{equation*}
\mathrm{Hom}_{\mathscr{D}_{Y}}^{*,*}(L_{0}, L_{1}) = \mathrm{ker }\, Q/
\mathrm{im }\, Q,
\end{equation*}
is defined by Floer theory, which itself is modeled after Morse theory
approach to supersymmetric quantum mechanics, from the introduction. The
role of the Morse complex is taken by the Floer complex.

For branes equipped with a trivial local system, the Floer complex
%
%e\upshape  5.7 #&#
\begin{equation}
\label{FC}
\mathrm{CF}^{*,*}(L_{0}, L_{1}) = \bigoplus _{{\mathcal P} \in L_{0} \cap L_{1}}
{\mathbb C} {\mathcal P}
\end{equation}
is a graded vector space spanned by the intersection points of the two
Lagrangians, together with the action of a differential $Q$. The complex
is graded by the fermion number, which is the Maslov index, and the equivariant
gradings, thanks to the fact $W$ is multivalued.

The action of the differential on this space
\begin{equation*}
Q: \mathrm{CF}^{*,*}(L_{0}, L_{1}) \rightarrow \mathrm{CF}^{*+1,*}(L_{0}, L_{1})
\end{equation*}
is generated by instantons. In Floer theory, the coefficient of
${\mathcal P}'$ in $Q{\mathcal P}$ is obtained by ``counting'' holomorphic
strips in $Y$ with boundary on $L_0$ and $L_1$, interpolating from ${\mathcal P}$ to
${\mathcal P}'$, of Maslov index $1$ and equivariant degree 0. The cohomology
of the resulting complex is Floer cohomology.

%s5.3.4 #&#
\subsubsection{}
\label{197.sec5.3.4}

A simplification in the present case is that, just as branes have a description
in terms of the Riemann surface, so do their intersection points, as well
as the maps between them.

The theory that results is a generalization of Heegard--Floer theory, which
is associated to $^{L}{\mathfrak{g}} ={\mathfrak{g l}}_{1|1}$ and categorifies
the Alexander polynomial \cite{OS0b, OS0c}. Heegard--Floer
theory has target
$Y_{{\mathfrak{g l}}_{1|1}} = \mathrm{Sym}^{d}({\mathcal A})$, the symmetric product
of $d$ copies of ${\mathcal A}$. $Y_{{\mathfrak{g l}}_{1|1}}$ should be thought of as a configuration
space of fermions on the Riemann surface, as opposed to anyons for
$Y_{{\mathfrak{su}}_{2}}$; in particular, their top holomorphic forms differ.

While we so far assumed that $^{L}{\mathfrak{g}}$ is simply laced, the
$\mathscr{D}_{Y}$ has an extension to non-simply laced Lie algebras, as
well as ${\mathfrak{gl}}_{m|n}$ and ${\mathfrak{sp}}_{m|2n}$ Lie superalgebras,
described in \cite{A3, ALR}.

%s5.4 #&#
\subsection{Link invariants and equivariant mirror symmetry}
\label{197.sec5.4}

Mirror symmetry helps us understand exactly which questions we need to
ask to recover homological knot invariants from $Y$.

%s5.4.1 #&#
\subsubsection{}
\label{197.sec5.4.1}

Since $Y$ is the ordinary mirror of $X$, we should start by understanding
how to recover homological knot invariants from $X$, rather than
${\mathcal X}$. Every B-brane on ${\mathcal X}$ which is relevant for us
comes from a $B$-brane on $X$ via an exact functor
%
%e\upshape  5.8 #&#
\begin{equation}
\label{fwant}
f_{*}: {\mathscr{D}_{X}} \rightarrow {\mathscr{D}_{\mathcal X}},
\end{equation}
which interprets a sheaf ``downstairs'' on $X$ as a sheaf ``upstairs''
on ${\mathcal X}$. The functor $f_{*}$ is more precisely the right-derived
functor ${R}f_{*}$. Its adjoint
%
%e\upshape  5.9 #&#
\begin{equation}
\label{fwantB}
f^{*}: {\mathscr{D}_{\mathcal X}} \rightarrow {\mathscr{D}_{X}}
\end{equation}
is the left derived functor ${L}f^{*}$, and corresponds to tensoring with
the structure sheaf $\otimes {\mathcal O}_{X}$, and restricting. Adjointness
implies that, given any pair of branes on ${\mathcal X}$ that come from
$X$,
\begin{equation*}
{\mathcal F} = f_{*} F, \quad {\mathcal G} = f_{*}G,%
\end{equation*}
the Hom's upstairs, in $\mathscr{D}_{\mathcal X}$, agree with the Hom's
downstairs, in $\mathscr{D}_{X}$,
%
%e\upshape  5.10 #&#
\begin{equation}
\label{model}
\mathrm{Hom}_{{\mathscr D}_{{\mathcal X}}}({\mathcal F}, {\mathcal G}) = \mathrm{Hom}_{{
\mathscr D}_{{X}}}( f^{*}f_{*}F, {G}),
\end{equation}
after replacing $F$ with $f^{*}f_{*}F$. The functor $f^{*}f_{*}$ is not
identity on $\mathscr{D}_{X}$.

%s5.4.2 #&#
\subsubsection{}
%%LEAP%%%\label{197.sec5.4.2}
\label{s-updown}

The equivariant homological mirror symmetry relating
$\mathscr{D}_{\mathcal X}$ and $\mathscr{D}_{Y}$ is not an equivalence
of categories, but a correspondence of branes and Hom's which come from
a pair of adjoint functors $h_{*}$ and $h^{*}$, inherited from
$f_{*}$ and $f^{*}$ via the downstairs homological mirror symmetry:
\begin{center}
\begin{tikzcd}[row sep=3em,column sep=3em]
{\mathscr{D}_{\mathcal X}} \arrow[dr,shift right, swap, "h^{*}"]
\\
&{\mathscr{D}_{Y}}\arrow[lu,shift right,swap,"h_{*}"]
\end{tikzcd}
\end{center}
Alternatively, $h^{*}$ and $h_{*}$ come by composing the upstairs mirror
symmetry with a pair of functors
$k^{*}: \mathscr{D}_{\mathcal Y}\rightarrow \mathscr{D}_{Y}$ and
$k_{*}: \mathscr{D}_{Y}\rightarrow \mathscr{D}_{\mathcal Y}$, which are
mirror to $f^{*}$ and $f_{*}$. The functors $k^{*}, k_{*}$ come from Lagrangian
correspondences; their construction is described in joint work with McBreen,
Shende, and Zhou \cite{AMSZ}. The functor $k_{*}$ amounts to pairing a
brane downstairs, with a vanishing torus fiber over it; this is how Figure~\ref{f_Ycore} arises in our model example. The adjoint functors let us
recover answers to all interesting questions about ${\mathcal X}$ from
$Y$.

%s5.4.3 #&#
\subsubsection{}
\label{197.sec5.4.3}

For any simply laced Lie algebra $^{L}\mathfrak{g}$, the branes
${\mathcal U} \in \mathscr{D}_{\mathcal X}$ which serve as cups upstairs
are the structure sheaves of (products of) minuscule Grassmannians, as described in Section~\ref{s:Gras}. They come via the functor
$h_{*}$ from branes $I_{\mathcal U}\in \mathscr{D}_{Y}$ downstairs, on $Y$
\begin{equation*}
{\mathcal U} = h_{*} I_{\mathcal U},%
\end{equation*}
which are (products of) generalized intervals. A minuscule Grassmannian
$G/P_{i}$ is the $h_{*}$-image of a brane which is the configuration space
of colored points on an interval ending on a pair of punctures on
${\mathcal A}$ corresponding to representations $V_{i}$ and
$V_{i}^{*}$. The points are colored by simple positive roots in
$\mu _{i} + \mu _{i}^{*} = \sum _{a} d_{a,i} \, {^{L}e}_{a}$, and ordered
in the sequence by which, to obtain the lowest weight $\mu _{i}^{*}$ in
representation $V_{i}$, we subtract simple positive roots from the highest
weight $\mu _{i}$. Because $V_{i}$ and $V_{i}^{*}$ are minuscule representations,
the ordering and hence the brane $I_{\mathcal U}$ is unique, up to equivalence
and a choice of grading.
%
%f7 #&#
\begin{figure}[t!]
\begin{center}
\includegraphics{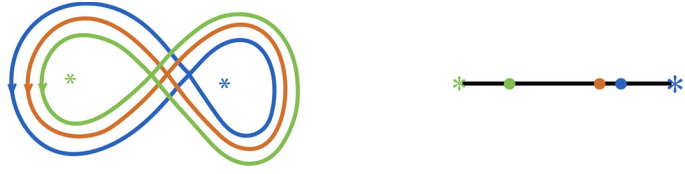}
\end{center}
%\vspace*{-6pt}
\caption{The cup and cap A-branes corresponding to the defining representation
of $^{L}{\mathfrak{g}} ={\mathfrak{su}}_{4}$, colored by its three simple
roots; they are equivariant mirror to a B-brane supported on a
${\mathbb P}^{4}$ as its structure sheaf.}%
\label{f-Eight}
\vspace*{-24pt}
\end{figure}
The ${\mathcal U}$ branes project back down as generalized figure-eight
branes; these are nested products of figure-eights, colored by simple
roots
\begin{equation*}
h^{*}{\mathcal U} = h^{*} h_{*} I_{\mathcal U} = E_{\mathcal U},%
\end{equation*}
and ordered analogously, as in Figure~\ref{f-Eight}. As objects of
$\mathscr{D}_{Y}$, these branes are best described iterated cones over
more elementary branes, mirror to stable basis branes \cite{ALR}. The cup
and cap branes all come with trivial local systems, for which the Floer
complexes are the familiar ones, given by~\eqref{FC}.

As an example, for $^{L}{\mathfrak{g}} = {\mathfrak{su}}_{2}$ the only
minuscule representation is the defining representation
$V_{i} = V_{\frac{1}{2}}$, which is self-conjugate. The cup brane
${\mathcal U}$ in ${\mathcal X}$ is a product of $d$ non-intersecting
${\mathbb P}^{1}$'s. It comes, as the image of $h_{*}$, from a brane
$I_{\mathcal U}$ in $Y$ which is a product of $d$ simple intervals, connecting
pairs of punctures that come together. The ${\mathcal U}$-brane projects
back down, via the $h^{*}$ functor, as a product of $d$ elementary figure-eight
branes. The branes are graded by Maslov and equivariant gradings, as described
in \cite{A2}.

%s5.4.4 #&#
\subsubsection{}
%%LEAP%%%\label{197.sec5.4.4}
\label{s_how}

In the description based on $Y$, both the branes, and the action of braiding
on them is geometric, so we can simply start with a link and a choice of
projection to the surface ${\mathcal A} = {\mathbb R} \times S^{1}$. A
link contained in a three ball in ${\mathbb R}^{2}\times S^{1}$ is equivalent
to the same link in ${\mathbb R}^{3}$, and projects to a contractible patch
on ${\mathcal A}$.

%f8 #&#
\begin{figure}[t!]
\vspace*{-6pt}
\begin{center}
\includegraphics{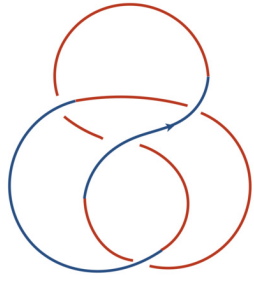}
\end{center}
\caption{A bicoloring of the left-handed trefoil.}%
\label{fT-1}
\vspace*{-18pt}
\end{figure}
To translate the link to a pair of A-branes, start by choosing bicoloring
of the link projection, such that each of its components has an equal number of red
and blue segments, and the red always underpass the blue. For a link component colored by a representation
$V_{i}$ of $^{L}{\mathfrak{g}}$, place a puncture colored by its highest
weight $\mu _{i}$ where a blue segment begins and its conjugate
$\mu _{i}^{*}$ where it ends; the orientation of the link component distinguishes
the two. The mirror Lagrangians $I_{\mathcal U}$ and
${\mathscr B} E_{\mathcal U}$ are obtained by replacing all the blue segments
by the interval branes, and all the red segments by figure-eight branes,
related by equivariant mirror symmetry to minuscule Grassmannian branes.
This data determines both $Y$ and the branes on it we need. The variant
of the second step, applicable for Lie superalgebras, is described in
\cite{ALR}.

Equivariant mirror symmetry predicts that a homological link invariant
is the space of morphisms
%
%e\upshape  5.11 #&#
\begin{equation}
\label{MP}
\mathrm{Hom}^{*,*}_{\mathscr{D}_{Y}}({\mathscr B} E_{\mathcal U},I_{\mathcal U})
= \bigoplus _{k\in {\mathbb Z}, {\vec d} \in {\mathbb Z}^{{rk}+1}}
\mathrm{Hom}_{\mathscr{D}_{Y}}\bigl({\mathscr B} E_{\mathcal U},I_{\mathcal U}[k]
\{ {\vec d}\}\bigr),
\end{equation}
the cohomology of the Floer complex of the two branes. In what follows,
will explain how to compute it.
%
%f9 #&#
\begin{figure}[h!]
\vspace*{-6pt}
\begin{center}
\includegraphics{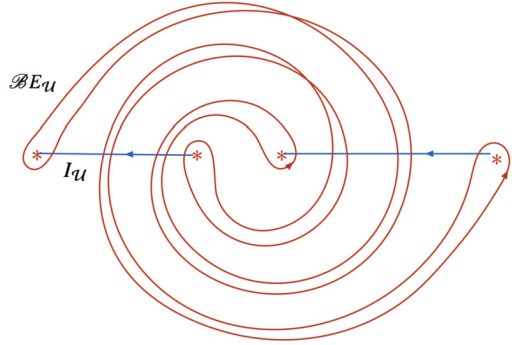}
\end{center}
\caption{The branes corresponding to the left-handed trefoil in
$^{L}{\mathfrak{g}} = {\mathfrak{su}}_{2}$. The knot was isotoped relative
to Figure \ref{fT-1}.}%
\label{f_T1}
\vspace*{-24pt}
\end{figure}
%

%s5.4.5 #&#
\subsubsection{}
\label{197.sec5.4.5}

To evaluate the Euler characteristic of the homology in~\eqref{MP}, one
simply counts intersection points of Lagrangians, keeping track
of gradings. For links in ${\mathbb R}^{3}$, the equivariant grading in~\eqref{MP} collapses to a ${\mathbb Z}$-grading. The Euler characteristic
becomes
%
%e\upshape  5.12 #&#
\begin{equation}
\label{EG}
\chi ({\mathscr B} E_{\mathcal U},I_{\mathcal U}) = \bigoplus _{{
\mathcal P} \in {\mathscr B} E_{\mathcal U} \cap I_{\mathcal U}} (-1)^{M({
\mathcal P})}{\mathfrak{q}}^{J({\mathcal P})},
\end{equation}
where $M({\mathcal P})$ and $J({\mathcal P})$ are the Maslov and
$c^{0}$-grading of the point ${\mathcal P}$; as in Heegard--Floer theory,
there are purely combinatorial formulas for them \cite{A3, ALR}. Mirror
symmetry implies that this is the
${U_{\mathfrak{q}}}(^{L}{\mathfrak{g}})$ invariant of the link in
${\mathbb R}^{3}$.

The fact that for $^{L}{\mathfrak{g}} = {\mathfrak{su}}_{2}$ the graded
count of intersection points in~\eqref{EG} reproduces the Jones polynomial
is a theorem of Bigelow \cite{bigelow}, building on the work of Lawrence
\cite{Lawrence1, Lawrence2, Lawrence3}. Bigelow also proved the statement
for $^{L}{\mathfrak{g}} = {\mathfrak{su}}_{N}$ with links colored by the
defining representation \cite{bigelowN}. The equivariant homological mirror
symmetry explains the origin of Bigelow's peculiar construction, and generalizes
it to other $U_{\mathfrak{q}}(^{L}{\mathfrak{g}})$ link invariants.\footnote{In \cite{M}, Bigelow's representation of the Jones polynomial, specialized to ${\mathfrak q}=1$, was related to the Euler characteristic of symplectic Khovanov homology of \cite{SS}.}

%s5.4.6 #&#
\subsubsection{}
\label{197.sec5.4.6}

The action of the differential $Q$ on the Floer complex, defined by counting
holomorphic maps from a disk ${D}$ to $Y$ with boundaries on the pair
of Lagrangians, should have a reformulation \cite{A2} in terms of counting
holomorphic curves embedded in ${D} \times {\mathcal A}$ with certain
properties, generalizing the cylindrical formulation of Heegard--Floer
theory due to Lipshitz \cite{L}. The curve must have a projection to
${D}$ as a $d=\sum _{a} d_{a}$-fold cover, with branching only between
components of one color, and a projection to ${\mathcal A}$ as a domain
with boundaries on one-dimensional Lagrangians of matching colors. In addition,
the potential $W$ must pull back to ${D}$ as a regular holomorphic function.
Computing the action of $Q$ in this framework reduces to solving a sequence
of well defined, but hard, problems in complex analysis in one dimension,
which are applications of the Riemannian mapping theorem, similar
to that in~\cite{OS0b}.

%s6 #&#
\section{Homological link invariants from A-branes}
\label{197.sec6}

To compute the link homology groups
%
%e\upshape  6.1 #&#
\begin{equation}
\label{MP2}
\mathrm{Hom}^{*,*}_{\mathscr{D}_{Y}}({\mathscr B} E_{\mathcal U},I_{\mathcal U}),
\end{equation}
we will make use of mirror symmetry which relates $X$ and $Y$ and is
the equivalence of categories
%
%e\upshape  6.2 #&#
\begin{equation}
\label{Downstairs}
\mathscr{D}_{X}\cong \mathscr{D}_{A}\cong \mathscr{D}_{Y},
\end{equation}
proven in \cite{ADZ}. A basic virtue of mirror symmetry is that it
sums up holomorphic curve counts. In our case, it solves all the disk-counting
problems required to find the action of the differential $Q$ on the Floer
complex underlying~\eqref{MP2}.

%s6.1 #&#
\subsection{The algebra of A-branes}
\label{197.sec6.1}

As in the simplest examples of homological mirror symmetry,
$\mathscr{D}_{X}$ and ${\mathscr{D}_{Y}}$ are both generated by a finite
set of branes.

%s6.1.1 #&#
\subsubsection{}
\label{197.sec6.1.1}

From perspective of $Y$, the generating set of branes come from products
of real line Lagrangians on ${\mathcal A}$, colored by
$d = \sum _{a} d_{a}$ simple roots. The brane is unchanged by reordering
a pair of its neighboring Lagrangian components, provided they are colored by roots which
are not linked in the Dynkin diagram
$\langle ^{L}e_{a} ,^{L}e_{b} \rangle =0$. It is also unchanged by passing
a component colored by $^{L}e_{a}$ across a puncture colored by a weight
$\mu _{i}$ with $\langle ^{L}e_{a} ,\mu _{i}\rangle =0$.

There is a generating brane
\begin{equation*}
T_{\mathcal C} = T_{i_{1}} \times \dots \times T_{i_{d}} \in
\mathscr{D}_{Y},%
\end{equation*}
for every inequivalent ordering of $d$ colored real lines on the cylinder.
Their direct sum
\begin{equation*}
T=\bigoplus _{\mathcal C} T_{\mathcal C} \in \mathscr{D}_{Y},%
\end{equation*}
is the generator of $\mathscr{D}_{Y}$ which is mirror to the tilting vector
bundle on $X$, which generates ${\mathscr{D}_{X}}$. This generalizes the
simplest example of mirror symmetry from Section~\ref{s-mirrorm}.
As before, we will be denoting branes on $X$ and on $Y$ related by homological
mirror symmetry by the same letter.
%\vspace*{-6pt}
%s6.1.2 #&#
\subsubsection{}
%%LEAP%%%\label{197.sec6.1.2}
\label{s-B}

A well known phenomenon in mirror symmetry is that it may introduce Lagrangians
with an extra structure of a local system, a nontrivial flat $U(1)$ bundle.
The mirror of a structure sheaf of a generic point, in our model example
of mirror symmetry from Section~\ref{s-mirrorm}, is a Lagrangian of this
sort.

Here, we find a generalization of this \cite{ADZ}. The pair of adjoint functors $h_*$ and $h^*$ that relate $\mathscr{D}_{Y}$ with its equivariant mirror $\mathscr{D}_{\mathcal{X}}$ equip each $T$-brane
with a vector bundle or, more precisely, with a local system of modules
for a graded algebra ${\mathcal B}$. The algebra is a product
${\mathcal B} = \bigotimes _{a=1}^{rk} {\mathcal B}_{{d_{a}}}$, where
${\mathcal B}_{{d}}$ has a representation as the quotient of the algebra
of polynomials in $d$ variables $z_{1}, \ldots , z_{d}$ which sets their
symmetric functions to zero. The $z$'s have equivariant
${\mathfrak{q}}$-degree equal to one.

As a consequence, the downstairs theory is not simply the Fukaya category of $Y$, but a Fukaya category
with coefficients \cite{ADZ}: Floer complexes assign to each intersection
point ${\mathcal P} \in L_{0} \cap L_{1}$ a vector space
$\hom_{{\mathcal B}}(\Lambda _{0}|_{\mathcal P}, \Lambda _{1}|_{\mathcal P})$ of homomorphisms of ${\mathcal B}$-modules
$\Lambda _{0,1}$ which $L_{0, 1}$ are equipped with. The
cup and cap branes $ E_{\mathcal U}$ and $I_{\mathcal U}$ come with trivial
modules for ${\mathcal B}$. The $T_{\mathcal C}$ branes correspond to modules
for ${\mathcal B}$ equal to ${\mathcal B}$ itself.
%\vspace*{-6pt}
%s6.1.3 #&#
\subsubsection{}
\label{197.sec6.1.3}

Since the $T_{\mathcal C}$-branes are noncompact, defining the Hom's between
them requires care. The Hom's
\begin{equation*}
\mathrm{Hom}_{\mathscr{D}_{Y}}\bigl(T_{\mathcal C}, T_{\mathcal C}' [k]\{\vec d\}\bigr) =
\mathrm{HF}\bigl(T_{\mathcal C}^{\zeta}, T_{\mathcal C}' [k]\{\vec d\}\bigr)
\end{equation*}
are defined through a perturbation of ${T}_{\mathcal C}$ which induces
wrapping near infinities of ${\mathcal A}$, as in Figure~\ref{f_cyl}, and
other examples of wrapped Fukaya categories.

The Floer cohomology groups $\mathrm{HF}$ are cohomology groups of the Floer complex
whose generators are intersection points of the $T_{\mathcal C}$ branes,
with coefficients in ${\mathcal B}$. The generators all have homological
degree zero, so the Floer differential is trivial, and
%
%e\upshape  6.3 #&#
\begin{equation}
\label{zero}
\mathrm{Hom}_{\mathscr{D}_{Y}}\bigl(T_{\mathcal C}, T_{{\mathcal C}'}[k]\{\vec d\}\bigr) =0,
\quad \textup{for all}\ k\neq 0\ \textup{and all}\ {\vec d\,}.
\end{equation}
The Floer product on $\mathscr{D}_{Y}$ makes
\begin{equation*}
A = \mathrm{Hom}^{*}_{\mathscr{D}_{Y}}(T, T) = \bigoplus _{{\mathcal C}, {
\mathcal C}'} \bigoplus _{{\vec d} \in {\mathbb Z}^{{rk}+1}}\mathrm{Hom}_{\mathscr{D}_{X}}\bigl(T_{\mathcal C}, T_{{\mathcal C}'}\{\vec d\}\bigr)
\end{equation*}
into an algebra, which is an ordinary associative algebra, graded only
by equivariant degrees.

%s6.1.4 #&#
\subsubsection{}
\label{197.sec6.1.4}

The vanishing in~\eqref{zero} mirrors the tilting property of $T$ viewed
as the generator of~$\mathscr{D}_{X}$. The tilting vector bundle
$T\in \mathscr{D}_{X}$ is inherited from the Bezrukavnikov--Kaledin tilting
bundle ${\mathcal T}$ on ${\mathcal X,}$
\begin{equation*}
{\mathcal T}=\bigoplus _{\mathcal C} {\mathcal T}_{\mathcal C} \in
\mathscr{D}_{\mathcal X},
\end{equation*}
from Section~\ref{s-BK}, as the image of the $f^{*}$ functor, which is tensoring
with the structure sheaf of $X$ and restriction,
$f^{*} {\mathcal T} = T \in \mathscr{D}_{X}$.
The endomorphism of the upstairs tilting generator ${\mathcal T,}$
\begin{equation*}
{\mathscr A}= \mathrm{Hom}^{*}_{\mathscr{D}_{\mathcal X}}({\mathcal T}, {
\mathcal T}),
\end{equation*}
is the cylindrical KLRW algebra.

Since ${\mathcal T}$ is a vector bundle on ${\mathcal X}$, the center of
${\mathscr A}$ is the algebra of holomorphic functions on
${\mathcal X}$. The downstairs algebra is a quotient of the upstairs one,
by a two-sided ideal
%
%e\upshape  6.4 #&#
\begin{equation}
\label{quotient}
A = {\mathscr A}/{\mathscr I}.
\end{equation}
The ideal ${\mathscr I}$ is generated by holomorphic functions that vanish on
the core $X$.

%s6.1.6 #&#
\subsubsection{}
\label{197.sec6.1.6}

The cKLRW algebra ${\mathscr A}$ is defined as an algebra of colored strands
on a cylinder, decorated with dots, where composition is represented by
stacking cylinders and rescaling \cite{W2}. The local algebra relations
are those of the ordinary KLRW algebra from \cite{webster}. Placing the theory on the cylinder, it gets additional gradings by the winding
number of strands of a given color, corresponding to the equivariant $\Lambda $-action
on ${\mathcal X}$.

The elements of the algebra $A= {\mathscr A}/{\mathscr I}$ have a geometric
interpretation by recalling the Floer complex
$\mathrm{CF}^{*}(T_{\mathcal C}, T_{{\mathcal C}'})$ is fundamentally generated
by paths rather the intersection points. The $S^{1}$ of the algebra cylinder
is the $S^{1}$ in the Riemann surface ${\mathcal A}$; its vertical direction
parameterizes the path. The geometric intersection points of the $T$-branes
on ${\mathcal A}$ correspond to strings in $A$. The flat bundle morphisms,
a copy of ${\mathcal B}$ for each geometric intersection point, dress the
strings by dots of the same color. The algebra ${\mathcal B}$ is the quotient,
of the subalgebra of ${\mathscr A}$ generated by dots, by the ideal
${\mathscr I}$ of their symmetric functions.

%s6.2 #&#
\subsection{The meaning of link homology}
\label{197.sec6.2}

Since $T=\bigoplus _{\mathcal C} T_{\mathcal C}$ generates
$\mathscr{D}_{Y}$, like every Lagrangian in $ \mathscr{D}_{Y}$, the
${\mathscr B} {E}_{\mathcal U}$ brane has a description as a complex
%
%e\upshape  6.5 #&#
\begin{equation}
\label{resolution}
{\mathscr B} {E}_{\mathcal U} \cong \cdots \xrightarrow{t_{1}} {
\mathscr B} { E}_{1}(T) \xrightarrow{t_{0}} {\mathscr B} { E}_{0}(T),
\end{equation}
every term of which is a direct sum of $T_{\mathcal C}$-branes. The complex is the projective resolution of the ${\mathscr B} {E}_{\mathcal U}$ brane. It
describes how to get the
${\mathscr B} {E}_{\mathcal U} \in \mathscr{D}_{Y}$ brane by starting
with the direct sum brane
%
%e\upshape  6.6 #&#
\begin{equation}
\label{aproB}
{\mathscr B} { E}(T)=\bigoplus _{k} {\mathscr B} { E}_{k}(T) [k],
\end{equation}
with a trivial differential, and taking iterated cones over elements
$t_{k}\in A$. This deforms the differential to
\begin{equation}\label{QA}
Q_{A} =\sum _{k} t_{k} \in A,
\end{equation}
which takes
\begin{equation*}
Q_{A}: {\mathscr B} {E}(T)\rightarrow {\mathscr B} {E}(T)[1],%
\end{equation*}
as a cohomological degree $1$ and equivariant degree $0$ operator, which squares to zero $Q_{A}^2=0$ in the algebra $A$.

%s6.2.1 #&#
\subsubsection{}
\label{197.sec6.2.1}

The category of A-branes $\mathscr{D}_{Y}$ has a second, Koszul dual set
of generators, which are the vanishing cycle branes
$I = \bigoplus _{\mathcal C} I_{\mathcal C}$ of \cite{A2}. The $I$- and the
$T$-branes are dual in the sense that the only nonvanishing morphisms from
the $T$- to the $I$-branes are
%
%e\upshape  6.7 #&#
\begin{equation}
\label{IT}
\mathrm{Hom}_{\mathscr{D}_{Y}}(T_{\mathcal C}, I_{{\mathcal C}'}) = {
\mathbb C} \delta _{{\mathcal C}, {\mathcal C}'}.
\end{equation}
The $I_{\mathcal C}$-branes and the $T{\mathcal C}$-branes are, respectively, the simple and the projective modules
of the algebra {$A$}.
\subsubsection{}
\label{197.sec6.2.1b}

Among the $I$-branes, we find the branes
$I_{\mathcal U}\in \mathscr{D}_{Y}$ which serve as cups. This is a wonderful
simplification because it implies that from the complex in~\eqref{resolution}, we get for free a complex of vector spaces:
%
%e\upshape  6.8 #&#
\begin{equation}
\label{resolutionhom}
0 \rightarrow \hom_{A}\bigl( {\mathscr B} { E}_{0}(T), I_{\mathcal U}\{
\vec{d}\}\bigr) \xrightarrow{t_{0}} \hom_{A}\bigl( {\mathscr B} { E}_{1}(T), I_{\mathcal U}\{\vec{d}\}\bigr) \xrightarrow{t_{1}} \cdots .
\end{equation}
The maps in the complex~\eqref{resolutionhom} are induced from the complex
in~\eqref{resolution}. The cohomologies of this complex are the link homologies
we are after,
%
%e\upshape  6.9 #&#
\begin{equation}
\label{final}
\mathrm{Hom}_{\mathscr{D}_{Y}}\bigl({\mathscr B} E_{\mathcal U},I_{\mathcal U}[k]
\{ {\vec d}\}\bigr) = H^{k}\bigl(\hom_{A}({\mathscr B}E_{\mathcal U}, I_{\mathcal U})\bigr).
\end{equation}
\subsubsection{}
\label{197.sec6.2.1b_}
We learn that link homology captures only a small part of the geometry
of ${\mathscr B} E_{\mathcal U}$, the braided cup brane, or more precisely, of the complex that resolves it. Because the
$T$-branes are dual to the $I$-branes by~\eqref{IT}, almost all terms in
the complex~\eqref{resolutionhom} vanish.
The cohomology~\eqref{final} of small complex that remains is the $U_{\mathfrak q}(^L{\mathfrak g})$ link homology.

%s6.2.2 #&#
\subsubsection{}
\label{197.sec6.2.2}

The complex~\eqref{resolutionhom} itself has a geometric interpretation
as the Floer complex,
\begin{equation*}
\mathrm{CF}^{*,*}( {\mathscr B} { E}_{\mathcal U}, I_{\mathcal U}).
\end{equation*}
Namely, the vector space at the $k$th term of the complex
\begin{equation*}
\hom_{A}\bigl( {\mathscr B} { E}_{k}(T), I_{\mathcal U}\{\vec{d}\}\bigr)
\end{equation*}
is identified, by construction described in section \ref{197.sec6.3}, with that spanned by the intersection points of the $\mathscr{B} E_{\mathcal U}$ brane and the
$I_{\mathcal U}$ brane, of cohomological degree $[k]$ and equivariant degree
$\{\vec{d}\}$.

The maps in the complex
\begin{equation*}
\cdots \xrightarrow{t_{k-1}} \hom_{A}\bigl( {\mathscr B} { E}_{k}(T), I_{\mathcal U}\{\vec{d}\}\bigr) \xrightarrow{t_{k}} \hom_{A}\bigl( {\mathscr B} { E}_{k+1}(T),
I_{\mathcal U}\{\vec{d}\}\bigr) \xrightarrow{t_{k+1}} \cdots
\end{equation*}
encode the action of the Floer differential. A priori, computing these
requires counting holomorphic disk instantons. In our case, mirror symmetry~\eqref{Downstairs} has summed them up.

%s6.3 #&#
\subsection{Projective resolutions from geometry}
\label{197.sec6.3}

The projective resolution in~\eqref{resolution} encodes all the
$U_{\mathfrak{q}}(^{L}{\mathfrak{g}})$ link homology, and more. Finding the resolution
requires solving two problems, both in general difficult. We will solve simultaneously \cite{ALR}.

%s6.3.1 #&#
\subsubsection{}
\label{197.sec6.3.1}

The first problem is to compute which module of the algebra $A$ the brane
${\mathscr B} E_{\mathcal U}$ gets mapped to by the Yoneda functor
\begin{equation*}
L \in \mathscr{D}_{Y}\rightarrow \mathrm{Hom}_{\mathscr{D}_{Y}}^{*,*}(T, L)
\in \mathscr{D}_{A}.
\end{equation*}
This functor, which is the source of the equivalence
$\mathscr{D}_{Y}\cong \mathscr{D}_{A}$, maps a brane
$L$ to a right module for $A$, on which the algebra acts as
\begin{equation*}
\mathrm{Hom}_{\mathscr{D}_{Y}}^{*,*}(T, L) \otimes \mathrm{Hom}_{\mathscr{D}_{Y}}^{*}(T,
T) \rightarrow \mathrm{Hom}_{\mathscr{D}_{Y}}^{*,*}(T, L).
\end{equation*}
Evaluating this action requires counting disk instantons.

%s6.3.2 #&#
\subsubsection{}
\label{197.sec6.3.2}

The second problem is to find the resolution of this module, as in~\eqref{resolution}. The Yoneda functor maps the $T_{\mathcal C}$ branes
to projective modules of the algebra $ A$, so the resolution in~\eqref{resolution} is a projective resolution of the $A$ module corresponding
to the ${\mathscr B} E_{\mathcal U}$ brane. This problem is known to be
solvable, however, only formally so, by infinite bar resolutions.

%s6.3.3 #&#
\subsubsection{}
\label{197.sec6.3.3}

In our setting, these two problems get solved together. Fortune smiles
since the ${\mathscr B}E_{\mathcal U} \in \mathscr{D}_{Y}$ branes are
products of $d$ one-dimensional Lagrangians on ${\mathcal A}$, for which
the complex resolving brane~\eqref{resolution} can be deduced explicitly,
from the geometry of the brane.

%s6.3.4 #&#
\subsubsection{}
\label{197.sec6.3.4}

Take a pair of branes $L'$ and $L''$ on $Y$ which are products of
$d$ one-dimensional Lagrangians on ${\mathcal A}$, chosen to coincide up
to one of their factors. Up to permutation, we can take
\begin{equation*}
L' = L_{1}\times L_{2} \times \dots \times L_{d}, \quad L'' = L''_{1}
\times L_{2} \times \dots \times L_{d}.%
\end{equation*}
If $L_{1}'$ and $L_{1}''$ (which are necessarily of the same color) intersect
over a point $p \in L_{1}'\cap L_{1}''$ of Maslov index zero, we get a
new one dimensional Lagrangian $L_{1}$ which is a cone over $p$,
\begin{equation*}
L_{1} = \mathrm{Cone}(p) = L_{1}' \xrightarrow{p} L''_{1},%
\end{equation*}
as well as a new $d$-dimensional Lagrangian $L$ on $Y$ given by
%
%e\upshape  6.10 #&#
\begin{equation}
\label{Lpp}
L =L_{1}\times L_{2} \times \dots \times L_{d}.
\end{equation}
The Lagrangian is a cone over the intersection point ${\mathcal P}$ of
$L'$ and $L''$ which is of the form
%
%e\upshape  6.11 #&#
\begin{equation}
\label{pD}
{\mathcal P} =(p, \mathrm{id}_{L_{2}}, \dots , \mathrm{id}_{L_{d}}) \in L' \cap L'',
\end{equation}
and which also has Maslov index zero, $L=\mathrm{Cone}({\mathcal P})$.

Conversely, any $L$ brane which is of the product form in~\eqref{Lpp} can
be written as a complex \cite{AurouxSM}
%
%e\upshape  6.12 #&#
\begin{equation}
\label{LC}
L\cong L' \xrightarrow{{\mathcal P}} L''
\end{equation}
with an explicit map ${\mathcal P}$ coming from a one-dimensional intersection
point in one of its factors, as in~\eqref{pD}.

%s6.3.5 #&#
\subsubsection{}
\label{197.sec6.3.5}

To find the projective resolution of the ${\mathscr B} E_{\mathcal U}$ brane in~\eqref{resolution}, start by isotoping the  brane, by stretching it straight along the cylinder.

Let the brane break at the two infinities of ${\mathcal A}$, to get the direct sum brane ${\mathscr B} E(T)$ in~\eqref{aproB}, on which the complex is based. To find the maps in the complex, record how the brane breaks, iterating the above construction, one one-dimensional intersection point at the time. Every intersection point of the form~\eqref{pD} translates to a specific element of the algebra $A$ and a specific map in the complex. The result is a product of $d$ one-dimensional complexes, which describes factors of ${\mathscr B} E_{\mathcal U}$, and captures almost
all
the terms in the differential $Q_{A}$. The remaining ones follow, up to quasi-isomorphisms, by
asking that the differential closes
$Q_{A}^{2}=0$ in the algebra $A$. The fact that not all terms in $Q_A$ are geometric is a general feature of $d>1$ theories.

In practice, it is convenient to first
break the brane one of the two infinities of~${\mathcal A}$, and only then on the other. The branes at the intermediate
stage are images, under the $h^{*}$ functor, of stable basis branes
\cite{MO, ese} on $\mathscr{D}_{\mathcal X}$. The stable basis branes play a similar role to that of
Verma modules in category ${\mathcal O}$. The detailed algorithm is given in \cite{ALR}.
%

%
%s6.3.6 #&#
\subsubsection{}
\label{197.sec6.3.6}

As an example, take the left-handed
trefoil and $^{L}{\mathfrak{g}} = {\mathfrak{su}}_{2}$, which leads to the brane configuration from Figure~\ref{f_T1}. For simplicity,
consider the reduced knot homology, where the unknot homology is set to be trivial. As in Heegard-Floer theory, this corresponds to
erasing a component from the ${\mathscr B}E_{\mathcal U}$ and the $I_{\mathcal U}$ branes, and leads to Figure~\ref{f_T2}. This also brings us back
to the setting of our running example, where $Y$ is the equivariant mirror
to ${\mathcal X}$, the resolution of the $A_{n-1}$ surface singularity, with $n=4$.

%f10 #&#
\begin{figure}[h!]
\begin{center}
\includegraphics{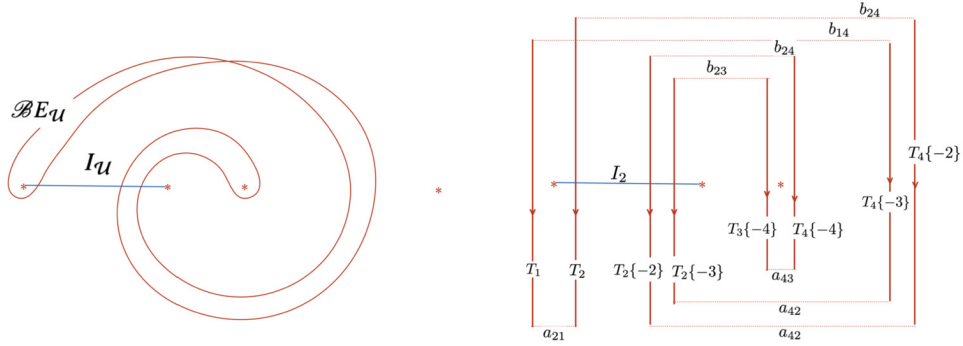}
\end{center}
\caption{Resolution of the ${\mathscr B}E_{\mathcal U}$ brane corresponding
to the reduced trefoil. The axis of the cylinder ${\mathcal A}$ is oriented
vertically here; the branes do not wind around the $S^1$.}%
\label{f_T2}
\end{figure}
The corresponding algebra $A=\bigoplus_{i,j=0}^{n-1}
\mathrm{Hom}_{{\mathscr D}_Y}^*(T_i,T_j)$ is the path algebra of an~affine
$A_{n-1}$ quiver, whose nodes correspond to $T_i$ branes.
The arrows $a_{i+1, i} \in \break\mathrm{Hom}_{{\mathscr D}_Y}(T_i,T_{i+1})$ and $b_{i, i+1} \in \mathrm{Hom}_{{\mathscr D}_Y}(T_{i+1},T_{i}\{1\})$ satisfy $a_{i,i-1}b_{i-1, i} =0= b_{i,i+1}a_{i+1, i} $, with $i$ defined modulo $n$. The $a$'s and $b$'s correspond to intersections of $T$-branes, near
one or the other infinity of ${\mathcal A}$; we have suppressed their $\Lambda$-equivariant degrees.

Isotope the ${\mathscr B}E_{\mathcal U}$ brane straight along the cylinder ${\mathcal A}$. Let it break
into $T$-branes, as in Figure~\ref{f_T2}, while recording how the brane breaks, one
connected sum at a time. Every connected sum of a pair of $T$-branes is
a cone over their intersection point, at one of the two infinities of
${\mathcal A}$, and a specific element of the algebra $A$. This leads to
the complex
%
%f11 #&#
\vspace*{9pt}
\begin{align*}
\includegraphics{197f11}
%\end{figure}
\end{align*}
\noindent which closes by the $A$-algebra relations.

The reduced homology of the trefoil is the cohomology of the
complex\break $\hom_{A}({\mathscr B}E^{\bullet}, I_{\mathcal U}\{d\})$
in~\eqref{resolutionhom}. The only non-zero contributions come from the $T_2$ brane, since the cup brane
$I_{\mathcal U} = I_{2}$ is dual to it. All
the maps evaluate to zero, as $I_{\mathcal U}$ brane is a simple module for $A$. As a consequence,
\begin{equation*}
\mathrm{Hom}_{\mathscr{D}_{Y}}\bigl({\mathscr B}E_{\mathcal U}, I_{\mathcal U}[k]
\{d\}\bigr) = H^{k}(\hom_{A}\bigl({\mathscr B}E^{\bullet}, I_{2}\{d\}\bigr)),%
\end{equation*}
equals to ${\mathbb Z}$ only for $(k,d) = (0,0), (2,-2), (3,-3)$, and vanishes
otherwise. Here, $k=M$ is the Maslov or cohomological degree and $d=J$ the Jones grading.
This is the reduced Khovanov homology of the left-handed trefoil, up to
regrading: Khovanov's $(i,j)$ gradings are related to $(M,J)$ by
$i = M+2J+i_{0}$ and $j = 2J+j_{0}$ where $i_{0}=0$,
$j_{0} = d+n_{+}-n_{-}$, where $n_{+}=0$, $n_{-}=3$ are the numbers of
positive and negative crossings, and $d=1$ is the dimension of~$Y$~\cite{A2}.

%s6.3.7 #&#
\subsubsection{}
\label{197.sec6.3.7}

The theory extends to non-simply-laced Lie algebras, and to Lie superalgebras
${\mathfrak{gl}}_{m|n}$ and ${\mathfrak{sp}}_{m|2n}$, as described in
\cite{ALR}.
The algebra $A$ corresponding to $^{L}{\mathfrak{g}}$ which is a Lie superalgebra, is not an ordinary associative algebra but a differential
graded algebra; the projective resolutions are then in terms of twisted complexes.
This section gives a method for solving the theory which is new even for $^{L}{\mathfrak{g}}={\mathfrak{gl}}_{1|1}$, corresponding to Heegard-Floer theory. The solution differs from that in \cite{OS}, in particular since our Heegard surface is ${\mathcal A}= {\mathbb R} \times S^{1}$, independent of the link.

\begin{ack}
This work grew out of earlier collaborations with Andrei Okounkov, which were indispensable. It includes
results obtained jointly with Ivan Danilenko, Elise LePage, Yixuan Li, Michael McBreen, Miroslav Rapcak, Vivek Shende,
and Peng Zhou. I am grateful to all of them for collaboration. I was supported by the NSF foundation grant PHY1820912,
by the Simons Investigator Award, and by the Berkeley Center for Theoretical
Physics.
\end{ack}
\end{document}